\documentclass[a4paper,reqno]{article}

\usepackage{amsmath}
\usepackage{amsfonts}
\usepackage{amssymb}
\usepackage{amsthm}
\usepackage[active]{srcltx}
\usepackage{color}
\usepackage{graphics, epsfig}
\usepackage{psfrag}
\usepackage[usenames,dvipsnames]{pstricks}

\allowdisplaybreaks

\newenvironment{ack}{\medskip{\it Acknowledgement.}}{}

\def\dashiint{\bint\kern-0.15cm\bint}
\renewcommand{\l}{\left}
\renewcommand{\r}{\right}

\newcommand{\pto}{(x_o,t_o)}

\begin{document}


\newtheorem{theorem}{Theorem}[section]
\newtheorem{proposition}{Proposition}[section]
\newtheorem{lemma}{Lemma}[section]
\newtheorem{corollary}{Corollary}[section]
\newtheorem{remark}{Remark}[section]
\newtheorem{definition}{Definition}[section]

\renewcommand{\thesection}{\arabic{section}}
\renewcommand{\theequation}{\thesection.\arabic{equation}}
\renewcommand{\thetheorem}{\thesection.\arabic{theorem}}
\numberwithin{equation}{section}
\numberwithin{theorem}{section}
\numberwithin{proposition}{section}
\numberwithin{lemma}{section}
\numberwithin{remark}{section}
\numberwithin{definition}{section}
\setcounter{secnumdepth}{3}



\newcommand{\cl}{\centerline}
\newcommand{\sms}{\smallskip}
\newcommand{\ms}{\medskip}
\newcommand{\bs}{\bigskip}
\newcommand{\noi}{\noindent}
\newcommand{\itl}[1]{\textit{#1}}
\newcommand{\blf}[1]{\textbf{#1}}
\newcommand{\dsty}{\displaystyle}
\newcommand{\txty}{\textstyle}
\newcommand{\ssty}{\scriptstyle}
\newcommand{\tty}{\texttt}


\newcommand\Par{\mathhexbox278\,}


\newcommand{\al}{\alpha}
\newcommand{\Al}{\Alpha}
\newcommand{\be}{\beta}
\newcommand{\Be}{\Beta}
\newcommand{\Gm}{\Gamma}
\newcommand{\gm}{\gamma}
\newcommand{\dl}{\delta}
\newcommand{\Dl}{\Delta}
\newcommand{\lm}{\lambda}
\newcommand{\Lm}{\Lambda}
\newcommand{\kp}{\kappa}
\newcommand{\varep}{\varepsilon}
\newcommand{\eps}{\epsilon}
\newcommand{\vp}{\varphi}
\newcommand{\sig}{\sigma}
\newcommand{\Sig}{\Sigma}
\newcommand{\om}{\omega}
\newcommand{\Om}{\Omega}
\newcommand{\uom}{\mbox{\boldmath$\omega$}}
\newcommand{\btau}{\mbox{\boldmath$\tau$}}
\newcommand{\bnu}{\mbox{\boldmath$\nu$}}
\newcommand{\up}{\upsilon}
\newcommand{\z}{\zeta}


\newcommand{\df}[1]{\buildrel\mbox{\small def}\over{#1}}
\newcommand{\op}[1]{\buildrel\mbox{\tiny o}\over{#1}}
\newcommand{\db}{\prime\prime}
\newcommand{\bsl}{\backslash}
\newcommand{\lb}{\lbrack\!\lbrack}
\newcommand{\rb}{\rbrack\!\rbrack}
\newcommand\la{\langle}
\newcommand\ra{\rangle}
\newcommand{\ev}{\equiv}
\newcommand{\nev}{\not\equiv}
\newcommand{\nn}{\mathbb{N}}
\newcommand{\qq}{\mathbb{Q}}
\newcommand{\zz}{\mathbb{Z}}
\newcommand{\rr}{\mathbb{R}}
\newcommand{\rn}{\rr^N}
\newcommand{\cc}{\mathbb{C}}
\newcommand{\id}{\mathbb{I}}
\newcommand{\bo}{\mathbb{O}}

\newcommand{\amsb}[1]{\mathbb{#1}}
\newcommand{\mcl}[1]{\mathcal{#1}}
\newcommand{\bl}[1]{\mathbf{#1}}
\newcommand{\ov}[1]{\overline{#1}}
\newcommand{\wt}[1]{\widetilde{#1}}
\newcommand{\wh}[1]{\widehat{#1}}

\newcommand{\lra}{\longrightarrow}
\newcommand{\LLR}{\Longleftrightarrow}
\newcommand{\LRA}{\Longrightarrow}
\newcommand{\LLA}{\Longleftarrow}


\newcommand{\bbox}{\vrule height.6em width.6em 
depth0em} 
\newcommand{\os}{\vbox{\hrule \hbox{\vrule 
height.6em depth0pt 
\hskip.6em \vrule height.6em depth0em}
\hrule}} 


\newcommand{\dvg}{\operatorname{div}}
\newcommand{\curl}{\operatorname{curl}}
\newcommand{\supp}{\operatorname{supp}}
\newcommand{\essup}{\operatornamewithlimits{ess\,sup}}
\newcommand{\essinf}{\operatornamewithlimits{ess\,inf}}
\newcommand{\essosc}{\operatornamewithlimits{ess\,osc}}
\newcommand{\osc}{\operatornamewithlimits{osc}}
\newcommand{\sign}{\operatorname{sign}}
\newcommand{\loc}{\operatorname{loc}}
\newcommand{\diam}{\operatorname{diam}}
\newcommand{\dist}{\operatorname{dist}}
\newcommand{\card}{\operatorname{card}}
\newcommand{\meas}{\operatorname{meas}}
\newcommand{\spn}{\operatorname{span}}
\newcommand{\dtm}{\operatorname{det}}
%


\newcommand{\overlim}{\mathop{\overline{\lim}}\limits}
\newcommand{\underlim}{\mathop{\underline{\lim}}\limits}
\newcommand{\ttop}[2]{\genfrac{}{}{0pt}{}{#1}{#2}}
\newcommand{\bcu}{\mathop{\txty{\bigcup}}\limits}
\newcommand{\bca}{\mathop{\txty{\bigcap}}\limits}
\newcommand{\bsu}{\mathop{\txty{\sum}}\limits}
\newcommand{\pro}{\mathop{\txty{\prod}}\limits}

\newcommand{\komjs}[1]{{\bf \blue \texttt{\symbol{'134}}#1\texttt{\symbol{'134}}}}


\newcommand{\pl}{\partial}
\newcommand{\ptt}{\frac{\pl}{\pl t}}
\newcommand{\ppx}{\frac\pl{\pl x}}
\newcommand{\dds}{\frac d{ds}}
\newcommand{\ddt}{\frac d{dt}}

\newcommand{\intl}{\int\limits}
\newcommand{\iintl}{\iint\limits}
\def\Xint#1{\mathchoice
    {\XXint\displaystyle\textstyle{#1}}%
    {\XXint\textstyle\scriptstyle{#1}}%
    {\XXint\scriptstyle\scriptscriptstyle{#1}}%
    {\XXint\scriptscriptstyle\scriptscriptstyle{#1}}%
    \!\int}
\def\XXint#1#2#3{\setbox0=\hbox{$#1{#2#3}{\int}$}
    \vcenter{\hbox{$#2#3$}}\kern-0.5\wd0}
\def\bint{\Xint-}
\def\dashint{\Xint{\raise4pt\hbox to7pt{\hrulefill}}}


\newcommand{\ovl}[3]{\int_{#1}^{#2}\kern-#3pt\raise4pt\hbox to7pt{\hrulefill}\ }

\newcommand{\ovll}[3]{\intl_{#1}^{#2}\kern-#3pt\raise4pt\hbox to7pt{\hrulefill}\ }

\newcommand{\tvl}[2]{\iint_{#1}\kern-#2pt\raise4pt\hbox to15pt{\hrulefill}\ }



\newcommand{\omt}{\Om_T}
\newcommand{\plo}{\partial\Omega}
\newcommand{\ovo}{\bar{\Om} }

%
\newcommand{\ci}[1]{C^\infty\!\left({#1}\right)}
\newcommand{\cio}[1]{C_o^\infty\!\left({#1}\right)}
\newcommand{\lloc}[1]{L_{\loc}\!\left({#1}\right)}
\newcommand{\xy}{|x-y|}


\newcommand{\intom}{\intl_{\Om}}
\newcommand{\intbo}{\intl_{\plo}}
\newcommand{\inom}{\int_{\Om}}
\newcommand{\inbo}{\int_{\plo}}
\newcommand{\intrn}{\intl_{\rn}}


\newcommand{\bye}{\end{document}}

\title{Sharp Regularity for Weak Solutions to the Porous Medium Equation}

\author{Ugo Gianazza\\
Dipartimento di Matematica ``F. Casorati"\\
Universit\`a di Pavia\\ 
via Ferrata 1, 27100 Pavia, Italy\\
email: {\tt gianazza@imati.cnr.it}
\and
Juhana Siljander\\
Department of Mathematics and Statistics\\ 
University of Jyv\"askyl\"a\\
P.O. Box 35, 40014 Jyv\"askyl\"a, Finland\\
email: {\tt juhana.siljander@jyu.fi}}

\date{\today}

\maketitle

\begin{abstract}
Let $u$ be a nonnegative, local, weak solution to the porous medium equation for $m\ge2$ in a space-time cylinder $\Om_T$. Fix a point $\pto\in\Om_T$: if the average 
\[
a\df=\dashint_{B_r(x_o)}u(x,t_o)\,dx>0,
\] 
then the quantity $|\nabla u^{m-1}|$ is locally bounded in a proper cylinder, whose center lies at time $t_o+a^{1-m}r^2$. {This implies that in the same cylinder the solution $u$ is H\"older continuous with exponent $\al=\frac1{m-1}$, which is known to be optimal}. Moreover, $u$ presents a sort of instantaneous regularisation, which we {quantify}.
\end{abstract}

\section{Introduction}
Let $\Omega$ be an open set in $\rn$ and for $T>0$ let $\Omega_T$ denote the cylindrical domain 
$\Omega\times(0,T]$, and $\Omega_{t_1,t_2}:=\Omega\times(t_1,t_2]$. 
Consider the well-known porous medium equation, namely the
quasilinear, degenerate parabolic partial 
differential equation
\begin{equation}\label{Eq:1:5:3}
\partial_t u-m\dvg (|u|^{m-1}\nabla u)=0,\quad m>1,\quad
\text{ weakly in }\> \Om_T.
\end{equation}
{The question of sharp regularity for equation~\eqref{Eq:1:5:3} has {been a long-standing open problem, which has divided expert opinions ever since the} famous counter-example by Aronson and Graveleau~\cite{graveleau,aron-grav}, described later, {proving that solutions do not admit the $C^{\frac1{m-1}}$--regularity of the Barenblatt fundamental solution. In this paper, we show that even in the case of the behavior represented by this counter-example, the solutions regularize instantaneously after a possible \textit{singular} point at the free boundary $\partial\{u>0\}$}, where the regularity fails. This is in the spirit of the classical time-lagged Harnack estimates for the heat equation, and with completely quantified local estimates.}

 {Moreover, the question of the sharp regularity is intimately related to the regularity of the free boundary. Caffarelli, V\'azquez and Wolanski~\cite{caff-vaz-wol} proved that under suitable conditions on the initial data, and after a suitable waiting time, the global solutions defined in $\rn$ admit the behavior we now obtain from purely interior arguments. Using the optimal gradient estimates, they continue to show the Lipschitz regularity of the free boundary, again after a waiting time depending on the initial data. 
 
 Contrary to their approach, we instead only consider locally defined quantities, and prove the optimal gradient estimates. While our aim has been to study the sharp regularity of the solutions, a significant question of its own, we also believe that similarly to~\cite{caff-vaz-wol}, {our result might be instrumental in concluding} the regularity problem with the free boundary for this highly nonlinear equation.}

{When proving regularity results}, at the heart of the matter lies the notion of solutions that one considers. We begin by giving the definition we will be working with. {For additional discussion of the existing literature and the meaning of our results, we refer to section~\ref{sec-nov-sig}.}
\begin{definition}\label{weak_solution}
A nonnegative function 
\begin{equation}\label{Eq:1:5:4}
\begin{array}{l}
u\in C_{\loc}\big(0,T;L^2_{\loc}(\Om)\big),\ \
u^{\frac{m+1}2}\in L^2_{\loc}\big(0,T; W^{1,2}_{\loc}(\Om)\big)\\
\end{array}
\end{equation}
is a local, weak sub(super)-solution to (\ref{Eq:1:5:3}) 
if for every compact set $K\subset \Om$ and every sub-interval 
$[t_1,t_2]\subset (0,T]$
\begin{equation}\label{Eq:1:5:5}
\begin{aligned}
\int_K u\vp\, dx\bigg|_{t_1}^{t_2}+\int_{t_1}^{t_2}\int_K
\big[&-u\,\partial_t \vp+ m u^{m-1}\nabla u\cdot \nabla\vp\big]dx\,dt\\
&\le(\ge)\,0
\end{aligned}
\end{equation}
for all nonnegative testing functions
\begin{equation}\label{Eq:1:5:6}
\vp\in W^{1,2}_{\loc}\big(0,T;L^2(K)\big)\cap 
L^2_{\loc}\big(0,T;W_o^{1,2}(K)\big).
\end{equation}
\end{definition}
\vskip.2truecm
This guarantees that all the integrals in (\ref{Eq:1:5:5}) 
are convergent. In \eqref{Eq:1:5:5} the symbol $\nabla u$ has to be understood in the sense of the following definition:
\begin{equation}\label{def:Du}
	\nabla u\df=\tfrac{2}{m+1}\mathbf 1_{\{ u>0\}} u^\frac{1-m}{2}\nabla u^\frac{m+1}{2}\, ,
\end{equation}
and for nonnegative functions
$$u^{m-1}\nabla u\df=\tfrac2{m+1} \mathbf 1_{\{ u>0\}} u^\frac{m-1}2\,\nabla u^\frac{m+1}2.$$

\begin{remark}
{\normalfont A number of different notions of solutions to the porous medium equation has been proposed in the literature. For a general overview and comparison among the different definitions, see for example \cite{vazquez-mono}. Here we stick to the one which has been widely used, when proving regularity estimates; see, for example, \cite{DBGV-mono,Ivanov-Mkrtychyan:1992}. The  requirements in \eqref{Eq:1:5:4} allow the testing of the  homogeneous equation only by the solution $u$ itself, and not by $u^m$, and lead to natural energy estimates for $u$ in terms of ${\nabla}u^\frac{m+1}{2}$, which in turn are instrumental in the proof of the Harnack inequality discussed in Theorem~\ref{Thm:3:15:1}. {We note, however, that the standard definition allowing testing with $u^m$ implies our assumptions on the solutions~\cite{ebmeyer}.}}
\end{remark}

For a local, globally bounded, weak solution $u:\Om_T\to\rr_+ {\cup \{0\}}$ to 
\eqref{Eq:1:5:3}, we are interested in the local boundedness of $|\nabla u^{m-1}|$  when $m\ge2$. The case $1<m<2$ has already been dealt with in \cite{Ivanov}, and we will not consider it here. Whenever such a bound on the gradient holds, it yields the sharp H\"older continuity exponent of $u$.

In~\cite{Ivanov} Ivanov proved the boundedness of $|\nabla u^{1+\kappa}|$, for some {large} $\kappa >0$. We will  use the same approach, but with some proper modifications, in order to prove the sharper result  $\kappa=m-2$.

\subsection{Notations}

{For $x_o\in \rn$ and $\rho > 0$, $B_{\rho}(x_o)$ denotes the ball of radius $\rho$, 
centered at $x_o$; when $x_o$ is the
origin of $\rn$, we simply write $B_{\rho}$.
We consider different kinds of cylinders; for $\theta > 0$ we define
\begin{equation*}
Q_{\rho}^{-}(\theta) = B_{\rho} \times (- \theta \rho^2,0], \quad Q_{\rho}^{+}
(\theta) = B_{\rho} \times (0,\theta \rho^2], \quad Q_{\rho}
(\theta) = B_{\rho} \times (-\theta\rho^2,\theta \rho^2],
\end{equation*}
where $\theta$ is a positive parameter that determines the length relative to $\rho^2$; 
for $\pto\in \rn \times \rr$
\begin{align*}
&\pto+Q_{\rho}^- (\theta) = B_{\rho} (x_o) \times (t_o-\theta \rho^2,t_o],\\
&\pto+Q_{\rho}^+ (\theta) = B_{\rho} (x_o) \times (t_o,t_o+\theta \rho^2],\\
&\pto+Q_{\rho} (\theta) = B_{\rho} (x_o) \times (t_o-\theta\rho^2,t_o+\theta \rho^2].
\end{align*}
Moreover, for $\rho,\theta>0$, and $\pto\in\rn\times\rr$ we let
\[
Q(\rho,\theta)=B_\rho\times(-\theta,0],\quad \pto+Q(\rho,\theta)=B_\rho(x_o)\times(t_o-\theta,t_o].
\]
Finally, we define
\begin{equation*}
U_{\rho}^{-}(\theta) = B_{\rho} \times (- \theta \rho,0], \quad 
U_{\rho}^{+}(\theta) = B_{\rho} \times (0,\theta \rho], \quad 
U_{\rho}(\theta) = B_{\rho} \times (-\theta\rho,\theta \rho],
\end{equation*}
where $\theta$ is a positive parameter that now determines the length relative to $\rho$;
for $\pto\in \rn \times \rr$
\begin{align*}
&\pto+U_{\rho}^- (\theta) = B_{\rho} (x_o) \times (t_o-\theta \rho,t_o]\\
&\pto+U_{\rho}^+ (\theta) = B_{\rho} (x_o) \times (t_o,t_o+\theta \rho]\\
&\pto+U_{\rho}(\theta) = B_{\rho} (x_o) \times (t_o-\theta\rho,t_o+\theta\rho). 
\end{align*}
When $\theta=1$, we write $Q_\rho(1)=Q_\rho$, and $U_\rho(1)=U_\rho$.

The parameters $\{m,N\}$ are the \emph{data}, and we say that a generic constant $\gm=\gm(m,N)$ depends upon the data, if it can be quantitatively determined a priori only in terms of the indicated parameters.

\subsection{The Intrinsic Scaling}
As it is always the case when dealing with degenerate equations such as \eqref{Eq:1:5:3},
an important step is to study the problem in the correct geometry: 
in our situation, it turns out that the natural cylinder for our purposes will be
\begin{equation}\label{Q_u}
Q_u \df=\pto+Q_{u\pto^{m-1}}(c u\pto^{1-m})=\pto+U_{u\pto^{m-1}}(c),
\end{equation}
where $c>0$ is a parameter that depends only on the data $\{m,N\}$.} In the following, we will talk of
\emph{intrinsic} cylinders, whenever their dimensions are determined by the value of the solution $u$; this can occur either through the value of $u$ at some point $\pto$, or through a proper integral average of $u$.
\subsection{The Main Result}
{
Let $Q=B \times I \Subset \Omega_T$ be a cylinder where $B$ is a ball and $I$ is a time interval. We choose another cylinder $Q^*=B^* \times I^* \Subset \Omega_T$ such that $Q \Subset Q^*$ and 
\begin{equation*}
\dist( Q^*, \partial \Omega_T), \dist(Q, \partial Q^*) \ge \frac14\dist(Q, \partial \Omega_T).
\end{equation*}
We also denote
\begin{equation}\label{def:do}
d_o \df =\frac12\dist(Q, \partial Q^*).
\end{equation}
Observe that $d_o$ is now controlled from below and above by $\dist(Q, \partial \Omega_T)$. Moreover, note that we consider the topological boundaries of $Q, Q^*$ and $\Omega_T$, and not just their parabolic boundaries.
}

\begin{theorem}\label{thm:main}
Let $u$ be a local, weak solution to \eqref{Eq:1:5:3} in $\Om_T$ with $m\ge 2$. Choose $\pto\in Q\Subset Q^*$.  For $r \in (0, \tfrac{d_o}2)$, assume that
\[
\dashint_{B_r(x_o)} u(x, t_o) \, dx>0,
\]
and define  
\begin{align*}
\theta_1 = \left[{\dsty \dashint_{B_r(x_o)} u(x, t_o) \, dx}\right]^{1-m}, \quad t_1 = t_o + \theta_1 r^2.
\end{align*}
Suppose that $\pto+Q^+_{4r}(\theta_1)\subset Q^*$.
Then there exist $C_o>1$, depending only on the data $\{m,N\}$, ${\dist(Q, \partial \Omega_T)}$, and the supremum of $u$ in $Q^*$, and $\bar\al>0$, depending only on the data $\{m,N\}$, such that 
\begin{align*}
\essup_{(x_o,t_1) + {Q}_{\frac r2}(\theta_1)}|\nabla u(x,t)^{m-1}| \le    C_o{\theta_1^{\bar\al}\l[\frac {d_o}r\r]^{2}}.
\end{align*} 
\end{theorem}
\noi Roughly speaking, for any point $\pto\in\Om_T$, provided a time-lag which depends only on the local behavior of $u$ is considered, $\nabla u^{m-1}$ is bounded in a proper cylinder. As a consequence, we have the following.
\begin{corollary}\label{Cor:hoelder}
Under the same assumptions of Theorem~\ref{thm:main}, for any compact set ${\mathcal K}\subset(x_o,t_1) + {Q}_{\frac r2}(\theta_1)$ the H\"older continuity exponent of $u$ in $\mathcal K$ is $\al=\frac1{m-1}$.
\end{corollary}
The gist of the previous results is that, given a point $\pto\in\Om_T$, we can conclude that the solution $u$ achieves its optimal modulus of continuity only after a finite, precisely quantified time lag. It turns out that we have a much better situation, and $u$ undergoes a sort of instantaneous regularisation, quantified by the next result.
\begin{theorem}\label{inst_reg_main}
Let $u$ be a a local, weak solution to \eqref{Eq:1:5:3} in $\Om_T$ with $m\ge 2$. Choose $\pto\in Q\Subset Q^*$. Let $\dsty {M = {1} + \sup_{Q^*}u}$. Then {for every $\gamma\ge 1$} at least one of the following holds:
\begin{enumerate}
\item {For every $r$ small enough we have}
\begin{equation}\label{m-1_holder-main}
\sup_{(x,t) \in \pto + Q_r({M^{1-m}})} |u(x,t)- u(x_o, t_o)| \le {\gm\, r^\frac1{m-1}}.
\end{equation}

\item There exist a sequence $\delta_n \to 0$ and a constant $C>0$, depending only on the data $\{m,N\}$, and ${\dist(Q, \partial \Omega_T)}$, such that 
\begin{equation}\label{eq_quan-main}
|\nabla u(x,t)^{m-1}| \le \frac{C}{{\gamma^{2(m-1)\bar\al}}\delta_n^{3\bar\al}},
\end{equation}
for every $\dsty(x,t) \in Q_n \df =  B_{\delta_n}(x_n) \times (t_o+\tfrac34\delta_n, t_o + \tfrac54\delta_n)$, where $\bar\al$ is quantity defined in Theorem~\ref{thm:main}. Consequently, by Corollary~\ref{Cor:hoelder}, we also have that $u \in C^{\frac1{m-1}}({\mathcal K})$, for every compact set ${\mathcal K}\subset Q_n$.
\end{enumerate}
\end{theorem}
\subsection{Novelty and Significance}\label{sec-nov-sig}
{When dealing with the porous medium equation, the optimal regularity of the pressure $p$ given by
\[
p\df=\frac m{m-1}u^{m-1},
\]
and consequently, by Darcy's law, of the velocity $\underline{v}$
\[
\underline{v}\df=-\frac m{m-1}\nabla u^{m-1}
\] 
is a problem that has attracted a lot of interest since the very beginning. When $N=1$, Aronson \cite{aron-SIAM} showed that $u^{m-1}$ is Lipschitz continuous with respect to $x$. In view of the Barenblatt fundamental solution
\begin{equation}\label{Barenblatt}
	{\mathcal B_m}(x,t)
	:=
	\left\{
	\begin{array}{cl}
	\frac1{t^{\frac{N}{\lm}}}
	\left[
	1-b\left(\frac{|x|}{t^{\frac1\lm}}\right)^2
	\right]_+^{\frac1{m-1}},&t>0\\[10pt]
	0 &t\le 0
	\end{array}
	\right.
\end{equation}
where 
\begin{equation}\label{def:lambda-b}
	\lm=N(m-1)+2
	\quad
	\mbox{and}
	\quad
	 b=b(m,N)=\frac{m-1}{2m\lm}\, ,
\end{equation} 
this is optimal. Still in the $1$-dimensional case, under a suitable monotonicity assumption on 
$(u^{m-1})_{xx}$, DiBenedetto \cite{DiBe79} proved that $u^{m-1}$ is Lipschitz continuous also with respect to time. Few years later, B\'enilan \cite{benilan}, and independently Aronson  and Caffarelli \cite{aron-caff}, proved the same Lipschitz continuity in $t$ without extra assumptions. This settles the problem when $N=1$. }

When $N>1$, the situation is much less clear. An overview of the results up to 1986 is in \cite{aron}. In particular, in \cite{benilan-1981} B\'enilan proved that $\nabla u^m$ is bounded if $(m-1)^2(N-1)<1$, showing that 
$$|\nabla u^m|^2\le c\frac ut,$$
where $c$ is a positive parameter, which depends only on the data $\{m,N\}$. 
In \cite{Moulay-Pierre}, given $Q=\rn\times(0,\infty)$, Moulay and Pierre showed that there exists $p>m-1$ such that $\nabla u^p\in L^{\infty}_{\loc}(Q)$ under the assumption 
$$\nabla u^r\in L^{N+2+\epsilon}_{\loc}(Q)$$ 
for any $\epsilon>0$ and some $r>0$, which is satisfied under suitable conditions on $\{m,N\}$. Their proof relies on classical iterative arguments of Moser type. 

Once one has a result of this kind, an estimate on the H\"older continuity exponent of $u$, both in space and in time, can be obtained by a result due to Kruzhkov \cite{kruzhkov}, later refined by Gilding \cite{gild76}. Under the point of view of the sharp H\"older continuity exponent, J\"ager and Lu \cite{Jaeger-Lu, Lu-Jaeger} have explicit examples. 

In order to prove boundedness of the gradient of a proper power of $u$, Ivanov and Mkrtychan \cite{Ivanov1992,Ivanov} worked the other way around, since they proved the boundedness of  $\nabla u^{1+\kappa}$, where $\kappa$ depends on $N$, $m$, and $\al$, the global H\"older continuity exponent of $u$. All these results are somewhat far from the {natural} regularity, namely that $\nabla u^{m-1}$ is locally bounded, as suggested by the Barenblatt fundamental solution \eqref{Barenblatt}: that this is indeed the case, provided that $t\ge T_o$, where $T_o$ depends only on the initial condition in $\rn$ was proved by Caffarelli, Vazquez  and Wolanski \cite{caff-vaz-wol}. 

Our results show that both the waiting time $T_o$, and that Lipschitz continuity of $u^{m-1}$ with respect to $x$ are purely local facts, independent of boundary and initial conditions, provided one works in a proper compact set $\mathcal K\Subset\Om_T$. 

It is well known that in general the result cannot hold up to the boundary: indeed, if one considers the domain
\[
\Om=\{x\in\rn:\ x_N>0\}
\]
and prescribes homogeneous boundary conditions on $\Gamma\df=\partial\Om\times(0,\infty)$, the function
\[ 
u_1(x,t)=x_N^{\frac 1m}
\]
is a stationary solution of \eqref{Eq:1:5:3} in $\Om_T$ and $\nabla u_1^{m-1}$ is not bounded in a neighborhood of $\Gamma$.

Moreover, it is clear from the proof that we work with {locally} bounded solutions: the \emph{local} boundedness of solutions to \eqref{Eq:1:5:3}
has been proved, for example, in \cite{andreucci,dahlberg-kenig}, but it is also well-known that solution can blow up in finite time, as shown by the weak solution $u_2:\rn\times(-\infty, T)\to\rr$, defined by
\[ 
u_2(x,t)=\left(\frac{m-1}{2\lm N}\right)^{\frac1{m-1}}\left(\frac{|x|^2}{T-t}\right)^{\frac1{m-1}},
\]
with $\lm$ as in \eqref{def:lambda-b}. 

Notice that our results do not contradict the famous counterexample discussed by Aronson and Graveleau \cite{graveleau,aron-grav}. They build a one-parameter family of radial self-similar solutions
$$v_c(x,t)=(T-t)^{2a-1}V_c\bigg(\frac{|x|}{(T-t)^a}\bigg)$$ for any $c>0$, where 
$$v={\frac m{m-1}}u^{m-1},$$ 
$$V_c(\xi)=\left\{
\begin{aligned}
0,&\qquad\text{if}\ \ 0\leq\xi\leq(\frac c{\gamma_m})^a\\ 
\text{positive,}&\qquad\text{if}\ \ \xi>(\frac c{\gamma_m})^a,
\end{aligned}
\right.$$ 
and $a=a_m\in(\frac12,1),\ \gamma_m\in(1,\infty)$ are constants. In a neighborhood of the interface, $v$ has a power-like behavior, whose exponent can only be numerically determined,
but such that $\nabla u^{m-1}$ is not bounded. Far from being a contradiction, this is exactly the second alternative considered in Theorem~\ref{inst_reg_main}. {Indeed, assuming the behavior of this counter-example at a point $\pto$, in addition to giving the optimal regularity for later times $t >t_o$, Theorem~\ref{inst_reg_main} also quantifies the blow-up of the optimal behavior around this point.}

{The structure of the paper is the following: in \S~\ref{section:preliminaries} we collect some preliminary results; \S~\ref{section:ivanov} is devoted to the regularisation of the function $u$, and to presenting the original Ivanov's argument: for the sake of completeness, we have included the proof; \S~\ref{S:Bound}--\ref{S:hoelder} deal with the regularised solution: first we prove Theorem~\ref{thm:main} for such a solution, and then we present the proof of Corollary~\ref{Cor:hoelder}; \S~\ref{S:limit} deals with the limit, as the regularisation parameter tends to zero, and completes the proof of Theorem~\ref{thm:main}. Finally, \S~\ref{S:Inst} concerns the \emph{Instantaneous Regularisation Theorem}. The last section takes into account some possible generalizations.}

{\begin{ack}
We acknowledge the warm hospitality of the Institut Mittag-Leffler, where this paper was started, during the program ``Evolutionary problems'' in the Fall 2013. Moreover, Juhana Siljander has been supported by Academy of Finland grant 259363 and a V\"ais\"al\"a foundation travel grant.
\end{ack}}
\section{Auxiliary results} \label{section:preliminaries}
\subsection{The Energy Estimate}
\begin{lemma}\label{energy_estimate}
Let $u \ge 0$ be a local, weak sub-solution to \eqref{Eq:1:5:3} in $\Om_T$. Then there exists a positive constant $\gamma=\gamma(m, N)$ such that for every cylinder $\pto+Q(\rho,\theta)\subset \Om_T$, every $k \ge 0$,  and every piecewise smooth, cutoff function $\varphi$ vanishing on the parabolic boundary $\partial_p[\pto+Q(\rho,\theta)]$ of the cylinder, we have
\begin{align*}
&\essup_{t_o-\theta \le \tau \le t_o} \int_{B_\rho(x_o)} (u-k)_+^{2} \varphi^2(x,\tau) \, dx + \iint_{\pto+Q(\rho,\theta)} u^{m-1}|\nabla (u-k)_+|^2 \, dxd\tau \\
&\le \gamma \iint_{\pto+Q(\rho,\theta)} u^{m-1}(u-k)_+^{2}|\nabla \varphi|^2 + (u-k)_+^{2}\l(\frac{\partial\varphi}{\partial t}\r)_+ \, dxd\tau.
\end{align*}
\end{lemma}
\vskip.2truecm
\begin{proof} For the proof we refer to~\cite[Proposition 6.1 of Chapter 3]{DBGV-mono}.
\end{proof}
\vskip.2truecm
\subsection{The Harnack Inequality}
Fix $\pto\in \Omega_T$ such that $u\pto>0$ and construct the cylinders 
\begin{equation*}
\pto+Q_{\rho} (\theta)\quad\text{ where }\quad
\theta=\l(\frac c{u\pto}\r)^{m-1},
\end{equation*} 
and ${c \ge 1}$ is a given positive constant. These cylinders 
are ``intrinsic'' to the solution, since their length 
is determined by the value of $u$ at $\pto$. In \cite[\S~11]{DBGV-acta} (see also \cite[Chapter~5]{DBGV-mono})
the following result has been proved.
\begin{theorem}\label{Thm:3:15:1} 
Let $u$ be a nonnegative, local, weak 
solution to equation \eqref{Eq:1:5:3} with $m\ge1$ in $\Om_T$. 
There exist positive constants $c$ and $\gamma_{o,s}$ depending 
only on the data $\{m,N\}$, such that for 
all cylinders $\pto+Q_{4\rho} (\theta)\subset\Omega_T$,
we have
\begin{equation}\label{Eq:3:15:2}
\gamma_{o,s}^{-1}\sup_{B_\rho(x_o)} u(\cdot, t_o-\theta\rho^2)\le 
u\pto\le\gamma_{o,s} \inf_{B_\rho(x_o)} u(\cdot,t_o+\theta\rho^2). 
\end{equation}
\end{theorem}
\begin{remark}
{\normalfont For Theorem~\ref{Thm:3:15:1} to hold, the large reference cylinder $\pto+Q_{4\rho}(\theta)$ needs to be contained in $\Om_T$; we can require a smaller reference cylinder, namely $\pto+Q_{k\rho(\theta)}$ with $1<k<4$, but this will affect the values of $c$ and $\gm_{o,s}$. In other words, the constants depend also on $d\df=\dist([\pto+Q_\rho(\theta)],{\partial}\Om_T)$: assuming $\pto+Q_{4\rho} (\theta)\subset\Omega_T$ amounts to bounding $d$ from below.}
\end{remark}

\subsection{The Weak Harnack Inequality}
\begin{theorem}{\bf \cite[Chapter~5]{DBGV-mono}, {\cite{Leht}}}\label{weak_harnack} 
Let $u$ be a nonnegative, local, weak 
super-solution to equation \eqref{Eq:1:5:3} with $m>1$ in $\Om_T$, $(x_o,t_o)\in\Om_T$, and assume that $B_{8\rho}(x_o)\subset \Om$. 
There exist positive constants $c$ and $\gamma_{o,w}$ depending 
only on the data $\{m,N\}$, such that 
we have
\begin{equation}\label{Eq:3:15:2bis}
\dashint_{B_\rho(x_o)}u(x,t_o) \, dx \le c\left(\frac{\rho^2}{T-t_o}\right)^{\frac1{m-1}}+\frac12\,\gamma_{o,w} \inf_{B_{4\rho}(x_o)} u(\cdot,t) 
\end{equation}
for all times
\[
t_o+\frac12\theta_1\rho^2\le t\le t_o+\frac32\theta_1\rho^2,
\]
where 
\[
\theta_1=\min\left\{c^{1-m}\frac{T-t_o}{\rho^2},\left(\dashint_{B_\rho(x_o)} u(x,t_o)\,dx\right)^{1-m}\right\}.
\]
\end{theorem}
\noi It is straightforward to see that if 
\[
t_o+(2c)^{m-1}\rho^2\left(\dashint_{B_\rho(x_o)}u(x,t_o)\,dx\right)^{1-m}\le T,
\]
then
\begin{equation}
\dashint_{B_\rho(x_o)}u(x,t_o) \, dx \le \gamma_{o,w} \inf_{B_{4\rho}(x_o)} u(\cdot,t) 
\end{equation}
for all
\[
t_o+\frac12\left(\dashint_{B_\rho(x_o)}u(x,t_o)\,dx\right)^{1-m}\rho^2\le t\le t_o+\frac32\left(\dashint_{B_\rho(x_o)}u(x,t_o)\,dx\right)^{1-m}\rho^2.
\]
Therefore, we have the following.
\begin{corollary}
Let $u$ be a nonnegative, local, weak 
super-solution to equation \eqref{Eq:1:5:3} with $m>1$ in $\Om_T$, $t_o\in(0,T)$ and assume that
\[
\dashint_{B_\rho(x_o)} u(x,t_o)\,dx>0.
\] 
There exist positive constants $c$ and $\gamma_{o,w}$ depending 
only on the data $\{m,N\}$, such that 
\begin{equation*}
\dashint_{B_\rho(x_o)}u(x,t_o) \, dx \le \gamma_{o,w} \inf_{B_{4\rho}(x_o)} u(\cdot,t) 
\end{equation*}
for all times
\[
t_o+\frac12\theta_1\rho^2\le t\le t_o+\frac32\theta_1\rho^2,\qquad \theta_1=\left(\dashint_{B_\rho(x_o)} u(x,t_o)\,dx\right)^{1-m}
\]
provided that $B_{8\rho}(x_o)\times(t_o,t_o+(2c)^{m-1}\theta_1\rho^2]\subset\Om_T$.
\end{corollary}
\subsection{Sup-Estimates}
\begin{lemma}\label{Lm:sup-DiBe}
Let {$m>1$ and} $u$ be a nonnegative, local, weak {sub}solution to equation \eqref{Eq:1:5:3} in $\Om_T$. There exists a constant $\gamma$ that depends only on the data $\{m,N\}$, such that for all cylinders $\pto+Q(\rho,\theta)\subset \Om_T$, and $\forall\,\sigma\in(0,1)$
\begin{align*}
\sup_{B_{\sigma\rho}(x_o)\times(t_o-\sigma\theta,t_o]} u&\le\gamma\frac{\theta/\rho^2}{(1-\sig)^{{N(m+1)+2}}}\left(\sup_{t_o-\theta<\tau<t_o}\dashint_{B_\rho(x_o)}u(x,\tau)\,dx\right)^m\\
&+\gamma\left(\frac{\rho^2}{\theta}\right)^{\frac1{m-1}}.
\end{align*}
\end{lemma}
\noindent This is a consequence of the following result from \cite[Lemma~4]{dahlberg-kenig}.
\begin{lemma}\label{Lm:sup-DK}
Let $u$ be a smooth positive {sub}solution to equation \eqref{Eq:1:5:3} in $\dsty S^*=\{(x,t):\ B_2\times(-4,0]\}$. Then
\[
\sup_{B_{\frac12}\times(-\frac14,0]}u\le\gamma\left\{1+\sup_{-4<s<0}\int_{B_2}u(x,s)\,dx\right\}^\eta,
\]
where $\gamma$ and $\eta$ depend only on the data $\{m,N\}$.
\end{lemma}
However, it is not straightforward to obtain Lemma~\ref{Lm:sup-DiBe} from Lemma~\ref{Lm:sup-DK}, and we could not find an explicit reference for the former. Therefore, we provide a full proof of Lemma~\ref{Lm:sup-DiBe}, to keep the paper self-contained.
\begin{proof}
Without loss of generality, we can assume $(x_o,t_o)=(0,0)$. By Remark~3.1 of \cite{andreucci} with $\lm=1$ we have
\begin{align*}
\sup_{Q(\sig\rho,\sig\theta)}u\le&\gm\frac{\theta/{\rho^2}}{(1-\sig)^{N+2}}\dashiint_{Q(\rho,\theta)}u^m\,dxd\tau+\gamma\left(\frac{\rho^2}{\theta}\right)^{\frac1{m-1}}\\            
                                             \le&\gm\frac{\theta/{\rho^2}}{(1-\sig)^{N+2}}\left[\dashiint_{Q(\rho,\theta)}u^{m+1}\,dxd\tau\right]^{\frac m{m+1}}+\gamma\left(\frac{\rho^2}{\theta}\right)^{\frac1{m-1}}.
\end{align*}
Let $E\subset\rn$ be an open, bounded set and
\begin{equation}\label{Eq:embedding}
u\in L^\infty(0,T;L^2(E)),\qquad u^{\frac{m+1}2}\in L^2(0,T;W^{1,2}_o(E)).
\end{equation} 
If $N=2$, consider the Gagliardo-Nirenberg multiplicative embedding inequality, see for example, \cite[Chapter~1, Theorem~2.1]{DiBe93} with
\[
v=u^{\frac{m+1}2},\quad q=\frac{2(m+2)}{m+1},\quad s=\frac{m+3}{m+1},\quad p=2,\quad\al=\frac{m+1}{2(m+2)}. 
\] 
We have
\[
\begin{aligned}
\int_E u^{m+2}\,dx&\le \gm\left[\int_E |Du^{\frac{m+1}2}|^2\,dx\right]^{\frac12}\int_E u^{\frac{m+3}2}\,dx\\
&=\gm\left[\int_E |Du^{\frac{m+1}2}|^2\,dx\right]^{\frac12}\int_E u^{\frac12}u^{\frac{m+2}2}\,dx\\
&\le\gm\left[\int_E |Du^{\frac{m+1}2}|^2\,dx\right]^{\frac12}\left[\int_E u\,dx\right]^{\frac12}\left[\int_E u^{m+2}\,dx\right]^{\frac12},
\end{aligned}
\]
which yields
\[
\int_E u^{m+2}\,dx\le\gm\int_E |Du^{\frac{m+1}2}|^2\,dx\int_E u\,dx.
\]
Integrating over $(0,T)$ proves
\[
\iint_{E_T}u^{m+1+\frac2N}\,dxd\tau\le\gm\left(\sup_{0<\tau<T}\int_E u(\cdot,\tau)\,dx\right)^{\frac2N}\iint_{E_T}|Du^{\frac{m+1}2}|^2\,dxd\tau
\]
for $N=2$.
For $N\ge3$,  
\begin{align*}
\iint_{E_T}u^{m+1+\frac2N}\,dxd\tau&=\iint_{E_T}u^{m+1} u^{\frac2N}\,dxd\tau\\
                                                             &\le\int_0^T\left[\int_E(u^{\frac{m+1}2})^{\frac{2N}{N-2}}\,dx\right]^{\frac{N-2}N}\left[\int_E u\,dx\right]^{\frac2N}d\tau\\
                                                             &\le\gm\left(\sup_{0<\tau<T}\int_E u(\cdot,\tau)\,dx\right)^{\frac2N}\iint_{E_T}|Du^{\frac{m+1}2}|^2\,dxd\tau.
\end{align*}
Fix $\sig\in(\frac12,1)$, and consider the increasing sequences
\begin{align*}
&\rho_n=\sig\rho+(1-\sig)\rho\sum_{i=1}^n 2^{-i},\qquad\theta_n=\sig\theta+(1-\sig)\theta\sum_{i=1}^n 2^{-i},\\
&Q_o=Q(\sig\rho,\sig\theta),\quad Q_\infty=Q(\rho,\theta),\quad Q_n=Q(\rho_n,\theta_n).  
\end{align*}
Let $(x,t)\to\zeta_n$ be a nonnegative, piecewise smooth, cut-off function in $Q_{n+1}$, that equals one in $Q_n$, vanishes on the parabolic boundary of $Q_{n+1}$ and such that
\[
|D\zeta_n|\le\frac{2^{n+1}}{(1-\sig)\rho},\qquad0\le\zeta_{n,t}\le\frac{2^{n+1}}{(1-\sig)\theta}.
\]
The function $(u\z_n)$ satisfies \eqref{Eq:embedding} in $Q_{n+1}$. Therefore,
\begin{align*}
\iint_{Q_{n+1}}(u\z_n)^{m+1+\frac2N}\,dxd\tau\le&\gm\left(\sup_{-\theta_{n+1}<\tau<0}\int_{B_{\rho_{n+1}}} (u\z_n)(\cdot,\tau)\,dx\right)^{\frac2N}\\
&\iint_{Q_{n+1}}|D(u\z_n)^{\frac{m+1}2}|^2\,dxd\tau.
\end{align*}
The constant $\gm$ depends only upon the data $\{m,N\}$, and it is independent of $\rho$, $\theta$, $n$. By Lemma~\ref{energy_estimate} we obtain
\begin{align*}
\dashiint_{Q_{n+1}}&|D(u\z_n)^{\frac{m+1}2}|^2\,dxd\tau\\
\le&\gm\frac{2^{2n}}{(1-\sig)^2\rho^2}\dashiint_{Q_{n+1}}u^{m+1}\,dxd\tau+\gm\frac{2^n}{(1-\sig)\theta}\dashiint_{Q_{n+1}}u^2\,dxd\tau\\
\le&\gm\frac{2^{2n}}{(1-\sig)^2\theta}\left[\frac\theta{\rho^2}\dashiint_{Q_{n+1}}u^{m+1}\,dxd\tau+\left(\dashiint_{Q_{n+1}}u^{m+1}\,dxd\tau\right)^{\frac2{m+1}}\right]\\
\le&\gm\frac{2^{2n}}{(1-\sig)^2\theta}\left[\frac\theta{\rho^2}\dashiint_{Q_{n+1}}u^{m+1}\,dxd\tau+\left(\frac{\rho^2}\theta\right)^{\frac2{m-1}}\right],
\end{align*}
where we have used H\"older's and Young's inequality to estimate the last term. Without loss of generality, we can assume that
\[
\frac{\theta}{\rho^2}\dashiint_{Q_{n+1}} u^{m+1}\,dxd\tau>\left(\frac{\rho^2}\theta\right)^{\frac2{m-1}},
\]
otherwise the Lemma becomes trivial. Let us now set
\[
X_n=\dashiint_{Q_n}u^{m+1}\,dxd\tau.
\]
We obtain the recursive inequality
\begin{align*}
X_n&=\dashiint_{Q_n}u^{m+1}\,dxd\tau\le{\gamma}\dashiint_{Q_{n+1}}(u\z_n)^{m+1}\,dxd\tau\\
&\le{\gamma}\left[\dashiint_{Q_{n+1}}(u\z_n)^{m+1+\frac2N}\,dxd\tau\right]^{\frac{N(m+1)}{N(m+1)+2}}\\
      &\le\gm\left[\dashiint_{Q_{n+1}}|D(u\z_n)^{\frac{m+1}2}|^2 dxd\tau\left(\sup_{-\theta_{n+1}<\tau<0}\int_{B_{\rho_{n+1}}}u(x,\tau)\,dx\right)^{\frac2N}\right]^{\frac{N(m+1)}{N(m+1)+2}}\\
      &\le\gm\left[\frac{2^{2n}}{(1-\sig)^2}\dashiint_{Q_{n+1}}u^{m+1}\,dxd\tau\right]^{\frac{N(m+1)}{N(m+1)+2}}\\
      &\ \left[\sup_{-\theta_{n+1}<\tau<0}\dashint_{B_{\rho_{n+1}}}u(x,\tau)\,dx\right]^{\frac{2(m+1)}{N(m+1)+2}}\\
      &=\gm\frac{2^\frac{2N(m+1)n}{N(m+1)+2}}{(1-\sig)^{\frac{2N(m+1)}{N(m+1)+2}}}\left[\dashiint_{Q_{n+1}}u^{m+1}\,dxd\tau\right]^{\frac{N(m+1)}{N(m+1)+2}}\\
      &\ \left[\sup_{-\theta_{n+1}<\tau<0}\dashint_{B_{\rho_{n+1}}}u(x,\tau)\,dx\right]^{\frac{2(m+1)}{N(m+1)+2}},
\end{align*}
that is,
\[
X_n\le\gm\frac{2^\frac{2N(m+1)n}{N(m+1)+2}}{(1-\sig)^{\frac{2N(m+1)}{N(m+1)+2}}} X_{n+1}^{\frac{N(m+1)}{N(m+1)+2}}\left(\sup_{-\theta_{n+1}<\tau<0}\dashint_{B_{\rho_{n+1}}}u(x,\tau)\,dx\right)^{\frac{2(m+1)}{N(m+1)+2}}.
\]
By the interpolation Lemma~5.2 of \cite[Chapter~2]{DBGV-mono}, we conclude that there exists $\gm$, depending only upon the data $\{m,N\}$, such that
\[
\dashiint_{Q(\sig\rho,\sig\theta)}u^{m+1}\,dxd\tau\le\frac\gm{(1-\sig)^{N(m+1)}}\left[\sup_{-\theta<\tau<0}\,\dashint_{B_\rho}u(x,\tau)\,dx\right]^{m+1}.
\]
Therefore, we conclude
\[
\sup_{Q(\sig\rho,\sig\theta)}u\le\gm\frac{\theta/{\rho^2}}{(1-\sig)^{N(m+1)+2}}\left[\sup_{-\theta<\tau<0}\,\dashint_{B_\rho}u(x,\tau)\,dx\right]^m+\gm\left(\frac{\rho^2}\theta\right)^{\frac1{m-1}}.
\]
\end{proof}
\subsection{The H\"older Continuity}
Local, weak solutions $u$ to equation \eqref{Eq:1:5:3} 
are locally H\"older 
continuous in $\Omega_T$, and 
permit one to exhibit a quantitative 
H\"older modulus of continuity. 

{For a wide class of equations, whose prototype is the porous medium one, this result is due to DiBenedetto and Friedman \cite{DiBeFried}. H\"older continuity for solutions to the Cauchy problem for \eqref{Eq:1:5:3} was established by Caffarelli and Friedman \cite{Caffarelli-Friedman}; their approach relies on the special property (of such \emph{global} solutions)
\[
u_t\ge -\frac kt u,\qquad \frac1k=m-1+\frac2N.
\]
Continuity of solutions to degenerate parabolic equations $\displaystyle u_t=\Delta(|u|^{m-1} u)$ was proved by Caffarelli \& Evans \cite{Caffarelli-Evans}, but the modulus of continuity implicit in their proof is essentially of logarithmic kind. For more details, see also \cite{daska-kenig-mono}.}

Let 
\begin{equation*}
\Gamma=\partial \Omega_T-\overline \Omega\times\{T\}
\end{equation*} 
denote the parabolic 
boundary of $\Omega_T$, and for a compact set $\mathcal K\subset \Omega_T$ 
introduce the {\it intrinsic}, parabolic $m$-distance from 
$\mathcal K$ to $\Gamma$ by 
\begin{equation*}
m-\dist(\mathcal K;\Gamma)\,\df{=}\,
\inf_{\ttop{(x,t)\in \mathcal K}{(y,s)\in\Gamma}}
\Big(|x-y|+\|u\|_{\infty,\Om_T}^{\frac{m-1}2}|t-s|^{\frac12}\Big).
\end{equation*} 
\begin{theorem}\label{Thm:3:16:1} 
Let $u$ be a bounded, local, weak solution to equation 
\eqref{Eq:1:5:3} in $\Om_T$. 
Then $u$ is locally H\"older continuous 
in $\Omega_T$, and there exist constants $\gamma_1>1$ and $\alpha\in(0,1)$ 
that can be determined a priori only in terms of the 
data $\{m,N\}$, such that for every compact set 
$\mathcal K\subset \Omega_T$,
\begin{equation*}
|u(x_1,t_1)-u(x_2,t_2)|\le\gamma_1\|u\|_{\infty,\Omega_T}
\Big(\frac{|x_1-x_2|+\|u\|_{\infty,\Omega_T}^{\frac{m-1}2}|t_1
-t_2|^{\frac12}}{m-\dist(\mathcal K;\Gamma)}\Big)^\alpha
\end{equation*}
for every pair of points $(x_1,t_1)$, and $(x_2,t_2)\in \mathcal K$. 
\end{theorem}
\begin{proof} Fix a point in $\Omega_T$, which up to a translation 
we take to be the origin of $\rr^{N+1}$. For $\rho>0$ consider 
the cylinder 
\begin{equation*}
Q_{\epsilon}=B_{\rho}\times(-\rho^{2-(m-1)\epsilon},0],
\end{equation*}
where $\epsilon\in(0,1)$ is to be determined, and set
\begin{equation*}
\mu^+_o=\sup_{Q_{\epsilon}}u,\qquad \mu^-_o=\inf_{Q_{\epsilon}}u,\qquad 
\omega_o=\osc_{Q_{\epsilon}}u=\mu^+_o-\mu^-_o.
\end{equation*} 
With $\omega_o$ at hand, construct now the cylinder of intrinsic 
geometry 
\begin{equation*}
Q_o=B_{\rho}\times(-\omega_o^{1-m}\rho^2,0]=Q^-_\rho(\om_o^{1-m}).
\end{equation*}
If $\omega_o\ge \rho^{\epsilon}$, then $Q_o\subset Q_\epsilon$. 
Theorem~\ref{Thm:3:16:1} is then a consequence 
of the following (see \cite{DBGV-mono}):
\begin{proposition}\label{Prop:3:16:1}
There exist constants $\gamma>1$, and $\epsilon,\varepsilon,\delta\in(0,1)$, 
that can be quantitatively determined only in terms of the 
data $\{m,N\}$ and independent of $u$ and $\rho$, 
such that, if $\omega_o\ge \rho^{\epsilon}$, setting $\rho_o=\rho$ and 
\begin{equation*}
\omega_n=\delta\omega_{n-1},\quad \rho_n=\varepsilon\rho_{n-1},\quad 
Q_n=Q^-_{\rho_n}(\omega_n^{1-m}),\quad\text{ for }\> n=0,1,\dots
\end{equation*}
we have $Q_{n+1}\subset Q_n$, and 
\begin{equation*}
\osc_{Q_n}u\le\omega_n.
\end{equation*}
\end{proposition}
\noindent A consequence of the previous Proposition is:
\begin{lemma}\label{hoelder}
There exist constants $\gamma_1>1$ and $\alpha\in(0,1)$, that can be determined a priori only in terms of the data $\{m,N\}$, such that for all cylinders $Q_r^-(\omega^{1-m})$ with 
$0<r\le\rho$,
we have
\begin{equation}\label{Eq:2:6}
\osc_{Q_r^-(\omega^{1-m})} u\le\gamma_1\omega\l(\frac r\rho\r)^\alpha,
\end{equation}
where $\omega$ is such that 
\begin{equation}\label{starting}
\osc_{Q_\rho^-(\omega^{1-m})} u \le \omega.
\end{equation}
\end{lemma}
\end{proof}
\subsection{The Comparison Principle}
The comparison principle for solutions to the porous medium equation is due to \cite{aronson-crandall-peletier,dahlberg-kenig,wu-zhao-li}.  For a full account of this result, see also \cite{vazquez-mono}. Here we give the statement, which we will use in the following.

Let $\Om\subset\rn$ be a smooth, bounded domain, and consider the following boundary value problems
\begin{equation}\label{Eq:5:7:7}
\left\{
\begin{aligned}
&u_1\in C(0,T; L^{2}(\Om))\cap W^{1,1}(0,T; L^1(\Om))\\
&u_1^{\frac{m+1}2}\in  L^2(0,T; W^{1,2}(\Om))\\
&\partial_t u_1-m\dvg(u_1^{m-1}\nabla u_1)=0\quad\text{ weakly in }\> \Om_T\\
&{u_1^{\frac{m+1}2} (\cdot,t)\bigm|_{\partial\Om}=g_1(\cdot,t)\in L^2\big(0,T; W^{\frac12,2}(\partial\Om)\big)}\\
&u_1(\cdot,0)=u_{o,1}\in L^{2}(\Om),
\end{aligned}
\right.
\end{equation}
\begin{equation}\label{Eq:5:7:8}
\left\{
\begin{aligned}
&u_2\in C(0,T; L^{2}(\Om))\cap W^{1,1}(0,T; L^1(\Om))\\
&u_2^{\frac{m+1}2}\in  L^2(0,T; W^{1,2}(\Om))\\
&\partial_t u_2-m\dvg(u_2^{m-1}\nabla u_2)=0\quad\text{ weakly in }\> \Om_T\\
&{u_2^{\frac{m+1}2} (\cdot,t)\bigm|_{\partial\Om}=g_2(\cdot,t)\in L^2\big(0,T; W^{\frac12,2}(\partial\Om)\big)}\\
&u_2(\cdot,0)=u_{o,2}\in L^{2}(\Om).
\end{aligned}
\right.
\end{equation}
with $g_i\ge0$, and $u_{o,i}\ge0$ in the proper sense. 
With these specifications, the Dirichlet data $g_i(\cdot,t)$ on $\partial \Om$ 
for almost every $t\in(0,T)$, are taken in the sense of the traces of functions 
in $W^{1,2}(\Om)$ and the initial data $u_{o,i}$ are taken in the sense of 
continuous functions in $t$ with values in $L^{2}(\Om)$.  
\begin{proposition}\label{prop:comparison}
Assume that
\begin{itemize}
\item $g_1$ and $g_2$ are continuous on $\overline{S_T}$, and $g_1\le g_2$,
\item $u_{o,1}$ and $u_{o,2}$ are continuous in $\overline\Om$, and $u_{o,1}\le u_{o,2}$,
\end{itemize}
where $S_T$ denotes the lateral boundary of $\Om_T$. 
Then, $u_1\le u_2$.
\end{proposition}
\subsection{A Stability Result}
{In addition to the comparison principle, which gives the local uniqueness of solutions, we also have the following stability result.
\begin{proposition}\label{prop:stability}
Under the same assumptions of Proposition~\ref{prop:comparison}, suppose for non-negative $g_1$ and $u_{1,o}$ and for $\eps>0$ that
\begin{align*}
&g_2(x,t)=g_1(x,t)+\eps,\qquad\text{ for all }\ (x,t)\in S_T\\
&u_{o,2}(x)=u_{o,1}(x)+\eps,\qquad\text{ for all }\ x\in\Om.
\end{align*}
Denoting $u \df=u_1$ and $u_\eps \df=u_2$, there exists a subsequence $u_{\eps_k}$ such that
\begin{equation}\label{stability}
u_{\eps_k} \to u
\end{equation}
uniformly in $C^o(\overline{\Om_T})$.
\end{proposition}
}
We postpone the proof of this stability proposition to \S~\ref{S:limit} and proceed with the proof of the main result. For a similar statement, see also \cite{abdulla-00}.
\begin{remark}\label{stab:remark}{\normalfont 
Proposition~\ref{prop:stability} guarantees that all the intrinsic cylinders we are going to build in the next sections, relying on a properly regularised solution $u_\epsilon$, are close both to one another and to the analogous one defined in terms of the original solution $u$, {provided $\eps$ is small enough. In particular, even if not explicitly mentioned, in the sequel we will only consider such values of $\eps$ that belong to the subsequence $\{\eps_k\}$ given by Proposition~\ref{prop:stability}.} }
\end{remark}
\section{Local Gradient Estimates for Smooth Solutions}\label{section:ivanov}
{Denote
\begin{equation}\label{sigma_o}
\sigma_o \df= \frac{d_o}{8M^{m-1}} \in (0,1),
\end{equation}
where 
\begin{equation}\label{sigma_o-bis}
{M = 1 + \sup_{Q^*}u},
\end{equation} 
and $d_o$ is the quantity introduced in \eqref{def:do}. }
\begin{remark}\label{Rmk:4:0}
{\normalfont Since the quantity $d_o$ can be arbitrarily small, there is no loss of generality in assuming $\sigma_o\in(0,1)$, as done in \eqref{sigma_o}.}
\end{remark}
{Consider $\pto\in Q$ and choose $\rho>0$ such that}
\begin{equation}\label{rho_bounds}
\sigma_o [u\pto+{1}]^{m-1} \le \rho\le  d_o.
\end{equation}
Let
\[
U \df=\pto+U_{\rho}, \quad \rho>0.
\]
Observe that the choice of $\rho$ guarantees that $\dsty U\Subset Q^*$. Moreover, $\sig_o$ takes into account the size of $u$ in a neighborhood of $\pto$, and the distance of $\pto$ from the parabolic boundary of $Q^*$. 

As it is well-known, when dealing with the regularity of $\nabla u$, the heart of the matter lies in the study of the free boundary of the solution,
that is, the set where $u$ vanishes. Now, a difficulty appears at points where the value of the function $u$ is zero, since at such points the intrinsic cylinder shrinks at its upper vertex. 

In order to avoid the problems with such a degeneracy, we study a \emph{regularised} 
Cauchy-Dirichlet problem. Namely, if $u$ is a local, weak solution to  \eqref{Eq:1:5:3} with $m\ge2$ in $\Om_T$, we let $u_\epsilon$ be the weak solution to  
\begin{equation}\label{regularised_equation}
\begin{split}
\begin{cases}
\partial_t u_\eps&= m\dvg (u_\eps^{m-1}\nabla u_\eps), \quad \text{in} \quad Q^*\\
u_\eps&= u+ \eps, \quad \text{on} \quad \partial_p Q^*,
\end{cases}
\end{split}
\end{equation}
where $\partial_p Q^*$ is the parabolic boundary of $Q^*$, and $\eps$ is the positive parameter
we have just fixed. 

Since $u \ge 0$, by Proposition~\ref{prop:comparison} we have that $u_\eps \ge \eps$ everywhere in 
$\overline{Q^*}$. Therefore, in $\overline{Q^*}$ equation \eqref{regularised_equation} can be seen as a particular 
instance of a linear parabolic equation with bounded 
and measurable coefficients. By known results 
(see, for example, \cite[Chapter~II]{LSU}), {since $u$ is bounded and locally H\"older continuous in $\Om_T$}, 
the solution to \eqref{regularised_equation} 
is bounded and globally H\"older continuous in $\overline{Q^*}$. Consequently, \eqref{regularised_equation}
can be regarded as a linear parabolic equation with 
bounded, and H\"older continuous coefficients. Again 
by classical theory (see \cite[Chapter~V]{LSU}), 
one can conclude that the solution is 
indeed locally smooth. 

Thus, we now avoid the problem with the degeneracy of the intrinsic cylinder and we may also differentiate the equation to obtain a priori estimates. In particular, the equation may now be written in  non-divergence form as
\begin{equation}\label{non-divergence_form}
\begin{split}
\partial_t u_\eps&= m\dvg(u_\eps^{m-1}\nabla u_\eps) \\
&=mu_\eps^{m-1}\Delta u_\eps+m(m-1)u_\eps^{m-2}|\nabla u_\eps|^2 \\
&\df ={\bf A}(u_\eps)\Delta u_\eps + {\bf A}'(u_\eps) |\nabla u_\eps|^2.
\end{split}
\end{equation}
We will prove suitable estimates for $u_\eps$, which are stable with respect to the
regularisation parameter $\eps$. Eventually, in \S~\ref{S:limit} we will pass to the limit, as $\eps\to0$.
\subsection{A First Estimate}
For the sake of notation, in what follows we will drop $\eps$ from the subscript of $u$, and denote $u_\eps$ with $u$. 
We recall an argument, which is originally due to Ivanov~\cite{Ivanov}.
\begin{lemma}\label{first_ivanov_lemma}
Let $u$ be a nonnegative, classical solution to the problem~\eqref{regularised_equation} in $U$. Let ${\bf A}(u):=mu^{m-1}$ and ${\bf B}(u):=2^6 N m^3u^{m-3}$. Let also $z=z(s)$ be a positive, nondecreasing $C^2$ function defined on the range of $u$ in $\overline{U}$. Consider the auxiliary function 
\[
\overline w \df= \frac{|\nabla u|^2\varphi}{z(u)},
\]
 where $\varphi$ is the cut-off function defined by
 \begin{equation}\label{cut_off}
\varphi(x,t)
=
\begin{cases}
(1-|x-x_o|^2\rho^{-2})^2(1-|t-t_o|^2\rho^{-2})^2 \quad \text{in} \quad U, \\
0 \quad \text{outside} \quad U.
\end{cases}
\end{equation}
Then, there exists a constant $C=C(m,N)>0$ such that for every $\theta >0$ we have that 
 either the function $\overline w$ attains its maximum in $\overline{U}$ at an interior point $(x^*, t^*) \in U$, and at this point we have that  
\begin{equation}\label{ivanov_estimate}
\mathcal J(\theta)\df=-z''(u)-\frac{z'(u)^2}{z(u)} -\l[\frac{{\bf B}(u)}{{\bf A}(u)}+\frac{C{\bf A}^{-1}(u)}{\rho \theta^2}+\frac C{\rho^2 \theta^2}\r] z(u) \le 0,
\end{equation}
or
\[
|\nabla u\pto| \le \l[\frac{z(u(x_o,t_o))}{z(u(x^*,t^*))}\r]^{\frac12}\theta.
\]
\end{lemma}
\begin{proof}
{Let $u_k=\partial_{x_k}u$; by applying} $u_k\partial_{x_k}$ to equation~\eqref{non-divergence_form}  we obtain the following for $v \df =|\nabla u|^2$. 
\begin{align*}
-\frac12\partial_t v +\frac12{\bf A}(u)\Delta v &={\bf A}(u) \sum_{k=1}^{N}\sum_{i=1}^{N} u_{ik}^2 -{\bf A}'(u)v\Delta u  \\
&\quad-{\bf A}''(u)v^2- 2{\bf A}'(u)\nabla u \cdot [\nabla^2 u \nabla u].
\end{align*}
Next we introduce the cut-off function {defined} in~\eqref{cut_off}. Denote $w \df =v\varphi$ and multiply the above inequality by $\varphi$ to obtain for $w$ that
\begin{align*}
&-\partial_t w + \frac{w\partial_t \varphi}{\varphi}+ {\bf A}(u)\l[\Delta w-\frac{w\Delta \varphi}{\varphi}-2\frac{\nabla w \cdot \nabla \varphi}{\varphi}+2\frac{w|\nabla \varphi|^2}{\varphi^2}\r] \\
&=2 {\bf A}(u)\sum_{k=1}^{N}\sum_{i=1}^{N}\varphi u_{ik}^2  -2{\bf A}'(u) \big [w\Delta u +2\varphi \nabla u \cdot [\nabla^2 u \nabla u] \big ] -2{\bf A}''(u)\frac{w^2}{\varphi} .
\end{align*}
Let $z$ and $\overline w$ be as in the statement of the lemma. We obtain
\begin{equation}\label{first_ivanov}
\begin{split}
&z(u)[-\partial_t \overline w+ {\bf A}(u)\Delta \overline  w]-z'(u)\overline w\partial_t u\\
&+{\bf A}(u)\big [2z'(u)\nabla u \cdot \nabla \overline w+ z''(u)\overline w |\nabla u|^2 + z'(u) \overline w\Delta u \big ] \\
&=2 {\bf A}(u)\varphi\sum_{k=1}^{N}\sum_{i=1}^{N} u_{ik}^2  +{\bf A}(u)\frac{w\Delta \varphi}{\varphi}   \\
&\quad+2{\bf A}(u)\l[\frac{[z'(u)\overline w \nabla u+z(u)\nabla \overline w] \cdot \nabla \varphi}{\varphi}-\frac{w|\nabla \varphi|^2}{\varphi^2}\r] \\
&\quad-2{\bf A}'(u)w\Delta u- 4{\bf A}'(u)\varphi \nabla u \cdot [\nabla^2 u \nabla u] \\
&\quad-2{\bf A}''(u)v w - \frac{w\partial_t \varphi}{\varphi}.\\
\end{split}
\end{equation}
We start estimating the different terms one by one. First of all, Young's inequality and Jensen's inequality give
\begin{align*}
2{\bf A}'(u)w\Delta u &= 2m(m-1)u^{m-2}w\Delta u \\
&\le \frac1{4N}mu^{m-1}\varphi\l[\sum_{i=1}^{N}u_{ii}\r]^2 + 4Nm(m-1)^2u^{m-3}\frac{w^2}{\varphi} \\
&\le \frac14{\bf A}(u)\varphi\sum_{i=1}^{N}u_{ii}^2 + 4Nm^3u^{m-3}z(u)v \overline w \\
&\le \frac14{\bf A}(u)\varphi\sum_{k=1}^{N}\sum_{i=1}^{N}u_{ik}^2 + \frac14{\bf B}(u)z(u)v \overline w,
\end{align*}
where ${\bf B}(u)=64Nm^3u^{m-3}$. Next, a repeated use of Cauchy-Schwarz inequality yields
\[
\nabla u \cdot [\nabla^2 u \nabla u] \le |\nabla u|^2\l(\sum_{k=1}^{N}\sum_{i=1}^{N}u_{ik}^2\r)^{1/2}=v\l(\sum_{k=1}^{N}\sum_{i=1}^{N}u_{ik}^2\r)^{1/2}.
\]
Thus, again by Young's inequality we obtain that
\begin{align*}
4{\bf A}'(u)\varphi \nabla u \cdot [\nabla^2 u \nabla u] &=4m(m-1)u^{m-2}\varphi \nabla u \cdot [\nabla^2 u \nabla u] \\
&\le \frac14 mu^{m-1}\varphi \sum_{k=1}^{N}\sum_{i=1}^{N}u_{ik}^2 + 16m(m-1)^2u^{m-3}\frac{w^2}\varphi \\
&\le \frac14 {\bf A}(u)\varphi \sum_{k=1}^{N}\sum_{i=1}^{N}u_{ik}^2 + \frac14{\bf B}(u)z(u)v \overline w.
\end{align*}
On the other hand, since $u$ is a solution to the regularised equation, similarly as above we obtain that
\begin{align*}
z'(u)\overline w \partial_t u - {\bf A}(u)z'(u) \overline w\Delta u&= z'(u) \overline w \Delta u^m - {\bf A}(u)z'(u) \overline w\Delta u\\
&=z'(u) \overline w \big[{\bf A}(u)\Delta u+{\bf A}'(u)|\nabla u|^2 \big] \\
&\ \ \ -{\bf A}(u)z'(u) \overline w\Delta u \\
&=m(m-1)u^{m-2}z'(u)v\overline w \\
&\ge-\frac{mu^{m-1}}2 \frac{z'(u)^2}{z(u)}v\overline w - \frac{m(m-1)^2}2u^{m-3}z(u)v\overline w \\
&\ge-\frac{{\bf A}(u)}{2}\frac{z'(u)^2}{z(u)}v\overline w - \frac{{\bf B}(u)}{4} z(u) v \overline w,
\end{align*}
where we used Young's inequality to estimate the right-hand side. Moreover, by Young's inequality,
\begin{align*}
2{\bf A}(u)\frac{z'(u)\overline w \nabla u \cdot \nabla \varphi}{\varphi} &\le 2{\bf A}(u)\frac{z'(u)\overline w |\nabla u| |\nabla \varphi|}{\varphi} \\
&\le\frac{{\bf A}(u)}{2}\frac{z'(u)^2}{z(u)}\, v\overline w+2\frac{{\bf A}(u)z(u)}{|\nabla u|^2}\l(\frac{|\nabla \varphi|}{\varphi}\r)^2v\overline w.
\end{align*}
Notice also that
\[
2{\bf A}''(u) \le \frac{{\bf B}(u)}4.
\]
Inserting the above estimates into~\eqref{first_ivanov} gives
\begin{equation}\label{second_ivanov}
\begin{split}
&z(u) \big [-\partial_t \overline w+ {\bf A}(u)\Delta \overline  w \big ]+2{\bf A}(u)\l[z'(u)\nabla u \cdot \nabla \overline w-\frac{z\nabla \overline w \cdot \nabla \varphi}{\varphi} \r] \\
&\ge -\l[z''(u)+\frac{z'(u)^2}{z(u)} +\frac{2z(u)}{|\nabla u|^2}\l(\frac{|\nabla \varphi|}{\varphi}\r)^2\r]{\bf A}(u)v \overline w\\
&\quad-{\bf B}(u)z(u)v \overline w  +2{\bf A}(u)\l[\frac{w\Delta \varphi}{2\varphi} -\frac{w|\nabla \varphi|^2}{\varphi^2}\r] - \frac{w\partial_t \varphi}{\varphi}. \\
\end{split}
\end{equation}
Next, we estimate the cut-off function. By direct calculations we have
\[
\frac{|\nabla \varphi|}{\varphi}\le\frac 4{\rho \varphi^{1/2}}, \quad\text{and}\quad \l(\frac{|\nabla \varphi|}{\varphi}\r)^2\le\frac {16}{\rho^2 \varphi},
\]
as well as
\[
\frac{|\Delta \varphi|}{\varphi}\le \frac{16N}{\rho^2\varphi} \quad \text{and} \quad\l|\frac{\partial_t \varphi}{\varphi}\r| \le \frac 4{\rho \varphi^{1/2}}.
\]
We use these estimates to control the terms containing $\varphi$ in~\eqref{second_ivanov}. We obtain 
\[
\frac{w\partial_t \varphi}{\varphi}  \le \frac{C{\bf A}^{-1}(u)z(u)}{\rho\varphi^{1/2} v}{\bf A}(u)v\overline w \le \frac{C{\bf A}^{-1}(u)z(u)}{\rho w}{\bf A}(u)v\overline w
\]
as well as
\begin{align*}
&2{\bf A}(u)\l[\frac{w\Delta \varphi}{2\varphi} -\frac{w|\nabla \varphi|^2}{\varphi^2}\r] \\
&\le \frac{2{\bf A}(u)z(u)}{|\nabla u|^2}\l[\frac{|\Delta \varphi|}{\varphi}+ \l(\frac{|\nabla \varphi|}{\varphi}\r)^2\r]v \overline w \le\frac{Cz(u)}{\rho^2w}{\bf A}(u)v \overline w.
\end{align*}
Finally, plugging these estimates into~\eqref{second_ivanov} yields
\begin{equation}\label{third_ivanov}
\begin{split}
&z(u) \big [-\partial_t \overline w+ {\bf A}(u)\Delta \overline  w \big ]+2{\bf A}(u)\l[z'(u)\nabla u \cdot \nabla \overline w-\frac{z\nabla \overline w \cdot \nabla \varphi}{\varphi} \r] \\
&\ge -\l[z''(u)+\frac{z'(u)^2}{z(u)} +\l[\frac{{\bf B}(u)}{{\bf A}(u)}+\frac{C{\bf A}^{-1}(u)}{\rho w}+\frac C{\rho^2 w}\r] z(u)\r]{\bf A}(u)v \overline w.
\end{split}
\end{equation}
Assume $(x^*, t^*)$ is the point of maximum for the function $\overline w$ in the set $\overline U$ and let $v_o:=v\pto, w_o:=w\pto$, $z_o:=z\pto$, as well as $w^*:=w(x^*, t^*)$ and $z^*:=z(x^*, t^*)$. Due to the choice of our cut-off function $\varphi$, $(x^*, t^*) \in U$, that is, it is an interior point of $U$. 

Given a $\theta >0$, we have the following alternative: either 
\begin{equation}\label{ivanov_alternative}
|\nabla u(x^*, t^*)| >0 \quad \text{and} \quad w^* > \theta^2,
\end{equation}
or one of these two inequalities does not hold. Assume first that both hold. In such a case, we have that
\[
\frac{C{\bf A}^{-1}(u)}{\rho w^*} \le \frac{C{\bf A}^{-1}(u)}{\rho \theta^2} \quad \text{and} \quad \frac{C}{\rho^2 w^*} \le \frac{C}{\rho^2 \theta^2}.
\]
Since $|\nabla u(x^*, t^*)|>0$, we have at this point that ${\bf A(u)} v \overline w >0$. 
As the left-hand side of~\eqref{third_ivanov} is nonpositive at $(x^*,t^*)$, we may use~\eqref{third_ivanov} to conclude~\eqref{ivanov_estimate}. 

On the other hand, if $|\nabla u(x^*, t^*)| =0$, we have directly that $\overline w(x^*,t^*)=0$, that is, 
$\overline w$ vanishes at its maximum point in $U$. Since this point lies in the interior of $U$, we have that $\varphi(x^*, t^*)>0$, which yields
\[
|\nabla u\pto|^2 = \frac{\overline w z}{\varphi} \le \frac{\overline w^* z}{\varphi}=0.
\]
Finally, if $w^* \le \theta^2$, we obtain directly that
\[
|\nabla u\pto|^2 = v_o=w_o \le \frac{z_o}{z^*}w^* \le  \frac{z_o}{z^*}\theta^2. 
\]
Altogether, we obtain in any case that either~\eqref{ivanov_estimate} holds, or
\[
|\nabla u\pto| \le \l(\frac{z_o}{z^*}\r)^{1/2}\theta.
\]
\end{proof}
\subsection{A Second Estimate}
In~\cite{Ivanov} Ivanov uses Lemma~\ref{first_ivanov_lemma} to conclude that $|\nabla u^{\kappa}|$ is bounded, where $\kappa \approx 1/\alpha$ with $\alpha \in (0,1)$ being the known H\"older exponent of $u$
given by \eqref{Eq:2:6} (see also \cite[Chapter~III]{DiBe93}, and \cite{DiBeFried}). After a clever choice of $z(u)$ in Lemma~\ref{first_ivanov_lemma}, the problem concerning the boundedness of the gradient actually reduces 
to proving an oscillation estimate for $u$. By relying on the known H\"older regularity of $u$, one may then conclude.

We will proceed in a similar way; however, instead of merely using the H\"older continuity of $u$, 
we aim at proving a new oscillation estimate in an intrinsic cylinder with a scaling similar to $Q_u$ introduced in~\eqref{Q_u}. Then, in the next sections, we will use this estimate to conclude the required regularity result.

Let $u_\eps$ be the classical solution to \eqref{regularised_equation}, 
and for $\delta\in(0,1]$ consider the cylinder
\[
Q_{u_\eps, \delta}\df=\pto+U_{\sigma_o[\delta u_\eps\pto]^{m-1}},
\]
where $\sig_o$ is the quantity defined in \eqref{sigma_o}. {By Proposition~\ref{prop:stability} (see also Remark~\ref{stab:remark}) we fix $\eps_o \in (0,1]$ small enough so that for every $\eps \le \eps_o$ we have}
\begin{equation}\label{Eq:3:10}
\begin{aligned}
\sig_o[\delta u_\eps\pto]^{m-1}&\le \sig_ou_\eps\pto^{m-1}\\
&\le \sig_o[u\pto+{1}]^{m-1}
\le\rho.
\end{aligned}
\end{equation}
Therefore,
\[
Q_{u_\eps, \delta}\subset U\Subset Q^*\Subset\Om_T.
\]
For the sake of notation, in what follows we will drop $\eps$ from the subscript of $u$, and denote $u_\eps$ with $u$. We now follow the argument of~\cite{Ivanov}, with some important modifications, and prove a second lemma. 

\begin{lemma}\label{Lm:3:2}
Consider $u$ in $\overline Q_{u, \delta}$. 
Then, there exist constants $\nu_o=\nu_o(m,N) \in (0,1)$ and ${C=C(m,N)>0}$ such that if
\begin{equation}\label{oscillation}
(1-\nu_o)u\pto \le u(x, t) \le (1+\nu_o)u\pto
\end{equation}
for every $(x,t) \in \overline Q_{u, \delta}$, we have that 
\[
|\nabla u\pto^{m-1}| \le {\frac{C}{\sigma_o\delta^{m-1}}}.
\]
In particular, the constant $C$ does not depend on the regularisation parameter $\eps$ of
\eqref{regularised_equation}.
\end{lemma}

\begin{proof}
We need to show that by assuming
\begin{equation}\label{oscillation_assumption}
(1-\nu)u\pto \le u(x, t) \le (1+\nu)u\pto=:m_2
\end{equation}
for some small $\nu>0$, quantitatively determined in terms of $N$ and $m$, we are able to show the boundedness of the gradient of $u^{m-1}$. The proof will determine how $\nu$ can be estimated below by a quantitatively determined $\nu_o=\nu_o(m,N)$.

In order to simplify the notation, we denote $u\pto$ with $u_o$, and let ${\bf A}(u)$ and ${\bf B}(u)$, 
$\overline w$ and $\mathcal J(\theta)$, be as in Lemma~\ref{first_ivanov_lemma}. 

We consider in $\overline Q_{u, \delta}$ the estimates we proved in Lemma~\ref{first_ivanov_lemma} for the general cylinder $U$. Suppose $\overline w$ attains its maximum in $\overline Q_{\delta, u}$ at $(x^*, t^*)$. By Lemma~\ref{first_ivanov_lemma} we know that there exists a constant $C=C(m,N)>0$ such that for every $\theta >0$ either
\begin{equation}\label{ivanov_estimate_duplicate}
\begin{split}
\mathcal J(\theta)=&-z''(u)-\frac{z'(u)^2}{z(u)} \\
& -\l[\frac{{\bf B}(u)}{{\bf A}(u)}+\frac{C{\bf A}^{-1}(u)}{\sigma_o[\delta u_o]^{m-1} \theta^2}+\frac C{\sigma_o^{2}[\delta u_o]^{2(m-1)} \theta^2}\r] z(u) \le 0
\end{split}
\end{equation}
holds in $(x^*, t^*)$, or
\[
|\nabla u\pto| \le \l(\frac{z_o}{z^*}\r)^{1/2}\theta.
\]
Observe that here we used $\rho = \sigma_o[\delta u_o]^{m-1}$ in the definition of $\mathcal J(\theta)$. This choice of $\rho$ might not satisfy the lower bound of~\eqref{rho_bounds}, but since~\eqref{ivanov_estimate} always implies the same estimate with smaller $\rho$, we in particular obtain~\eqref{ivanov_estimate_duplicate}.

We start by estimating $\mathcal J(\theta)$. For this purpose, we assume that~\eqref{oscillation_assumption} holds with $\nu=\nu_1 = 1/2$. With this assumption we have
\[
\frac{{\bf B}(u)}{{\bf A}(u)}= \frac{2^6 N m^3u^{m-3}}{mu^{m-1}}= 2^6 N m^2 u^{-2} \le 2^8 N m^2 u_o^{-2} \df =C_o(m,N)u_o^{-2}
\]
for every $(x,t) \in \overline Q_{u, \sigma}$. On the other hand, since $\delta, \sigma_o \le1$, we have
\begin{align*}
\frac{C{\bf A}^{-1}(u)}{\sigma_o[\delta u_o]^{m-1}\theta^2}+\frac C{\sigma_o^{2}[\delta u_o]^{2(m-1)} \theta^2} &= \frac{Cm^{-1}u^{1-m}}{\sigma_o[\delta u_o]^{m-1}\theta^2} +\frac C{\sigma_o^{2}[\delta u_o]^{2(m-1)} \theta^2} \\
&\le  \l[\frac{C2^{m-1}m^{-1}}{\sigma_o^{2}\delta^{2(m-1)}}+\frac{C}{\sigma_o^{2}\delta^{2(m-1)}}\r]\frac{1}{u_o^{2(m-1)}\theta^2} \\
&\df =\frac{C_1(m,N)}{(\delta u_o)^{2(m-1)}[\sigma_o\theta]^2}. 
\end{align*}
Notice that if \eqref{oscillation_assumption} holds with another $\nu^*<\frac12$, the previous estimates continue to hold.

Let $\omega\df=2\nu_2 u\pto$, where $0<\nu_2\le\frac12$ is to be chosen, and now assume that~\eqref{oscillation_assumption} holds with $\nu=\nu_2$. Our task is to show that $\nu_2$ can be chosen quantitatively. 
Now the aim is to use Lemma~\ref{first_ivanov_lemma} and for that purpose we choose
\[
z(u)=6\omega^2-(m_2-u)^2,
\]
where $m_2$ is as introduced in~\eqref{oscillation_assumption}.
Observe that by the assumption that~\eqref{oscillation_assumption} holds for $\nu=\nu_2$, 
we obtain that $m_2-u \le \omega$.
Now using~\eqref{ivanov_estimate_duplicate} and plugging in the choice of $z$ and the above estimates, we obtain
\begin{equation}\label{ivanov_contradiction}
2-\frac45-6C_o(m,N)\l(\frac\omega {u_o}\r)^2 -{6C_1(m,N)\frac{\omega^2}{(\delta u_o)^{2(m-1)}[\sigma_o\theta]^2} } \le\mathcal J(\theta) \le 0.
\end{equation}
for every $\theta>0$. Next, we choose $\nu_2=\nu_2(m,N)\le\frac12$ small enough, so that
\[
6C_o(m,N)\l(\frac\omega {u_o}\r)^2=24C_o(m,N)\nu_2^2 \le \frac13.
\] 
Finally, we choose 
\[
\theta = \frac{[48C_1(m,N)]^{1/2}{\nu_2\,}}{\sigma_o\delta^{m-1}u_o^{m-2}}\df=\frac{C_2(m,N)}{\sigma_o\delta^{m-1}u_o^{m-2}},
\]
to conclude that
\[
6C_1(m,N)\frac{\omega^2}{(\delta u_o)^{2(m-1)}[\sigma_o\theta]^2}  = \frac12.
\]
Plugging these into~\eqref{ivanov_contradiction} yields $\txty\frac{11}{30} \le 0$, and thus~\eqref{ivanov_estimate_duplicate} cannot hold. By Lemma~\ref{first_ivanov_lemma} we conclude that
\[
|\nabla u\pto| \le \l(\frac65\r)^{1/2}\theta
\]
which gives the claim, due to our choice of $\theta$, provided we let $\nu_o\df=\min\{\nu_1, \nu_2\}$. Notice that in the last estimate we used the fact that $5\omega^2 \le z(u) \le 6 \omega^2$.
\end{proof}
\section{The Local Bound for the Gradient of the Regularised Solution After a Waiting Time}\label{S:Bound}
The key point in Lemma~\ref{Lm:3:2} is the oscillation control on $u_\eps$ given by \eqref{oscillation}. We now proceed with showing how this control can be obtained after a quantifiably determined waiting time.
For the sake of notation, in what follows we will drop $\eps$ from the subscript of $u$, and denote $u_\eps$ with $u$.

Suppose $u$ is the positive, bounded solution to \eqref{regularised_equation} in $\overline{Q^*}$.
Let 
\[
\gamma_o\df=\max\{\gm_{o,s},\gm_{o,w}\}
\]
where $\gm_{o,s}$ and $\gm_{o,w}$ are the constant respectively in \eqref{Eq:3:15:2} and \eqref{Eq:3:15:2bis}. We have $\gm_o=\gamma_o(m,N, d_o)\ge 1$; similarly, let $\alpha_o=\alpha_o(m,N) \in (0,1)$ be the local H\"older continuity exponent and 
$\gamma_1=\gamma_1(m,N, d_o)\ge 2$ be the local H\"older continuity constant in $Q^*$, as given by Theorem~\ref{Thm:3:16:1} and Lemma~\ref{hoelder}. 

Recall the definition of $\sigma_o$ in~\eqref{sigma_o}--\eqref{sigma_o-bis} as
\[
\sigma_o= \frac{d_o}{8M^{m-1}} \in (0,1),\quad M = {1}+ \sup_{\overline {Q^*}} u.
\]
Let $u_o \df =u\pto$ with $\pto \in Q$. We have the following theorem.
\begin{theorem}\label{thm_quant}
Let $m\ge 2$, and suppose that $u$ is the $\eps$-regularised solution of~\eqref{regularised_equation} in $Q^* \Supset Q \ni \pto$ 
For $r \in (0, \tfrac{d_o}2)$, assume that
\[
\dashint_{B_r(x_o)} u(x, t_o) \, dx>0,
\]
and define  
\begin{align*}
\theta_1 = \left[{\dsty \dashint_{B_r(x_o)} u(x, t_o) \, dx}\right]^{1-m}, \quad t_1 = t_o + \theta_1 r^2.
\end{align*}
Suppose that $\pto+Q^+_{4r}(\theta_1)\subset Q^*$.
Then there exists a positive constant $C_o$, depending only on the data $\{m,N\}$, $d_o$, and $M$, but independent of the regularisation parameter $\eps$, such that 
\begin{equation}\label{eq:4.1}
\essup_{(x_o,t_1) + Q_{\frac r2}(\theta_1)}|\nabla u(x,t)^{m-1}| \le    C_o \theta_1^{\frac2{\alpha_o(m-1)} + 1}\l[\frac {d_o}r\r]^{2}.
\end{equation}
\end{theorem}

\begin{remark}
{\normalfont From the proof of Theorem~\ref{thm_quant} one can also check that the right-hand side of~\eqref{eq:4.1} is scaling invariant. Indeed, the constant $C_o$ here is of the form $$\gamma M^{\tfrac2{\alpha_o}+m-1}\sigma_o^2,$$
where $\gamma$ depends only on $\{m,N,\dist(Q, \partial\Omega_T)\}$. Now $\sigma_o$ is scaling invariant due to the intrinsic scaling corresponding to the equation, and $$ M^{\tfrac2{\alpha_o}+m-1}\theta_1^{\frac2{\alpha_o(m-1)} + 1}$$ is scaling invariant directly by definition.}
\end{remark}
Fix the radius $r \in \big(0, \tfrac {d_o}2\big)$, and let $\eta\in(0,1)$ be such that 
$$r =  \frac {\eta d_o}2.$$ 
With this notation we have the following lemma.
\begin{lemma}\label{first_nu}
Let $m\ge 2$, suppose that $u$ is the $\eps$-regularised solution of~\eqref{regularised_equation} in $Q^* \Supset Q \ni \pto$, and let $\nu_o$ be the constant of Lemma~\ref{Lm:3:2}. Under the same notation and assumptions of Theorem~\ref{thm_quant}, fix a point 
$(x^*, t^*) \in (x_o,t_1)+Q_{\frac r2}(\theta_1)$ and let $u_* \df =u(x^*,t^*)$.
 Then, for the constant
\begin{equation}\label{c_1}
c_1=\l[\frac{\nu_o}{\gamma_o\gamma_1 M}\dashint_{B_r(x_o)} u(x,t_o) \, dx\r]^{\frac2{\alpha_o}+m-1}\eta^2\,\sigma_o\in (0,1)
\end{equation}
we have that
\begin{equation}\label{u_lemma_bounds}
(1-\nu_o)u_* \le u(x,t) \le (1+\nu_o)u_*
\end{equation}
for every $(x,t) \in (x^*, t^*)+ U_{\rho}$ with $\rho \df=c_1\sigma_o u_*^{m-1}$.  
\end{lemma}

\begin{proof}
As we mentioned before, without loss of generality, we may assume $\sig_o\in(0,1)$. By the weak Harnack inequality of Theorem~\ref{weak_harnack} we have
\begin{equation}\label{wh}
\inf_{B_r(x_o)} u(\cdot,t) \ge \gamma_o^{-1} \dashint_{B_r(x_o)} u(x, t_o) \, dx
\end{equation}
for all $$t \in \left[t_o+\tfrac{1}2\theta_1r^2, t_o+\tfrac32\theta_1r^2\right].$$
Choose a point $\dsty(x^*,t^*) \in (x_o,t_1)+Q_{\frac r2}(\theta_1)$, and let
\begin{equation*}
\bar t \df = t^* + c_1 \sigma_o u_*^{m-1},
\end{equation*}
as well as
\begin{align*}
\bar \rho \df=  \l[\frac{\nu_o}{\gamma_o\gamma_1M}\dashint_{B_r(x_o)} u(x,t_o) \, dx\r]^{\frac1{\alpha_o}}\eta\sigma_ou_*^{m-1}.
\end{align*}
Observe that since
\begin{align*}
c_1  \sigma_ou_*^{m-1} &= \l[\frac{\nu_o}{\gamma_o\gamma_1 M}\dashint_{B_r(x_o)} u(x,t_o) \, dx\r]^{\frac2{\alpha_o}+m-1}\eta^2\sigma_o^2u_*^{m-1} \\
&\le \l[\frac{\eta d_o}{8M^{m-1}}\r]^2u_*^{m-1} \\
&\le \frac14 M^{1-m}\l[\frac{\eta d_o}{2}\r]^2,
\end{align*}
we have $\bar t \le t^*+\frac14\theta_1r^2$. We consider the intrinsic cylinder $(x^*,\bar  t) + Q^-_{\bar  \rho}(M^{1-m}) = B_{\bar  \rho}(x_o) \times (\bar t - M^{1-m}{\bar  \rho}^2, \bar t]$. Since $\sigma_o=\frac{d_o}{8M^{m-1}}$ and $\nu_o \in (0,1)$ we have that ${\bar  \rho} \le r/2$ and
\begin{equation*}
M^{1-m}{\bar  \rho}^2 \le {\tfrac14}\theta_1 r^2.
\end{equation*}
The H\"older continuity of $u$ stated in Lemma~\ref{hoelder}, together with the weak Harnack inequality~\eqref{wh}, now yields
\begin{equation}\label{cont}
\osc_{(x^*,\bar t) + Q^-_{\bar \rho}(M^{1-m})} u \le \gamma_1 M \left(\frac{{\bar  \rho}}{\sigma_o u_*^{m-1}}\right)^{\alpha_o} \le \frac{\nu_o}{\gamma_o}\dashint_{B_r(x_o)} u(x,t_o) \, dx \le \nu_o u_*.
\end{equation}
This implies
\begin{equation}\label{star_bounds}
(1-\nu_o)u_*\le u(x,t) \le (1+\nu_o)u_*
\end{equation}
for all $(x,t) \in (x^*,\bar t) + Q^-_{\bar \rho}(M^{1-m})$.

Due to the choice of $u_*$, we may apply the weak Harnack inequality~\eqref{wh} once more, to deduce for $\gamma_1 \ge 2$ that
\begin{align*}
M^{1-m}\bar \rho^2 &= M^{1-m}\l[\frac{\nu_o}{\gamma_o\gamma_1M}\dashint_{B_r(x_o)} u(x,t_o) \, dx\r]^{\frac2{\alpha_o}} [\eta\sigma_ou_*^{m-1}]^2 \\
&\ge 2\l[\frac{\nu_o}{\gamma_o\gamma_1M}\dashint_{B_r(x_o)} u(x,t_o) \, dx\r]^{\frac2{\alpha_o}+m-1} \eta^2\sigma_o^2u_*^{m-1} \\
&\ge 2c_1 \sigma_ou_*^{m-1}.
\end{align*}
This allows us to conclude that 
\begin{equation*}
(x^*, t^*) +  U_{ \rho} = B_{ \rho}(x^*) \times (t^*-\rho, t^*+\rho) \subset (x^*,\bar t) + Q^-_{\bar \rho}(M^{1-m}).
\end{equation*}
Therefore, we finally obtain 
\begin{align*}
(1-\nu_o)u(x^*, t^*) \le u(x,t) \le (1+\nu_o)u(x^*, t^*)
\end{align*}
for all $(x,t) \in (x^*, t^*) + U_{\rho}$.

\end{proof}

\begin{remark}
{\normalfont The solution $u$ is now bounded above and below in a neighbourhood of $(x^*, t^*)$, whose size depends only on $u_*$. Equivalently, we can say that $u$ is \emph{uniformly} small in a neighbourhood of $(x^*, t^*)$, whose dimensions are quantitatively known. In principle, this would not rule out the possibility for the gradient to be large, as $u$ might wildly oscillate close to $(x^*, t^*)$, independently of its size. However, since  the largeness of $\nabla u$ is controlled by the largeness of $u$ due to Lemma~\ref{Lm:3:2}, this implies that no large oscillation can actually occur.}
\end{remark}
We can now conclude.
\vskip.2truecm
\begin{proof}[Proof of Theorem~\ref{thm_quant}] By Lemma~\ref{first_nu}  we have bounded $u$ in terms of $u_*$ in a cylinder whose dimensions are proportional to $u_*^{m-1}$. Thus, in $(x^*, t^*) + U_{c_1\sigma_ou_*^{m-1}}$ we have the proper bounds required by Lemma~\ref{Lm:3:2}. Hence, there exists a positive constant $C$, which depends only on the data {$\{m,N\}$}, such that
\[
{|\nabla u(x^*, t^*)^{m-1}|\le \frac{C}{c_1\sigma_o}}.
\]
{Plugging in the choices of $\sigma_o$ and} $c_1$ allows us, due to the arbitrariness of $(x^*, t^*)\in (x_o, t_1) + {Q}_{\tfrac r2}(\theta_1)$, to conclude the proof.
\end{proof}
\section{The Optimal H\"older Continuity Exponent}\label{S:hoelder}
We now show how the local boundedness of $|\nabla u^{m-1}|$ can be used to prove a local $C^{\frac1{m-1}}$--H\"older estimate for $u$. The proof of Theorem~\ref{holder} below is originally due to Gilding~\cite{gild76}; here we modify it, in that we state it in intrinsic terms, in the spirit of Theorem~\ref{Thm:3:16:1}, although not exactly in the same way, as explained in Remark~\ref{Rmk:int:dist} below. 
Notice that by $C^{\frac1{m-1}}(Q)$, we denote what is usually referred to as $C^{\frac1{m-1},\frac1{2(m-1)}}(Q)$.
\begin{theorem}\label{holder}
Let $u$ be a nonnegative, local, weak solution to \eqref{Eq:1:5:3} in $\Omega_T$, with 
\begin{equation}\label{cont_bounds}
u(x,t) \le {\mathcal M}, \quad\text{and}\quad |\nabla u(x,t)^{m-1}| \le C_o
\end{equation}
for every $(x,t) \in  \Omega_T$. Then, for every $\mathcal K \Subset \Omega$, with $d_1 \df = \dist(\mathcal K, \partial \Omega)$, there exists a constant $\gamma=\gamma(m,N,C_o, d_1)>0$ such that
\begin{equation}\label{eq:hoelder}
|u(x, t) - u(y, s)| \le \gamma {\mathcal M}\l(|x-y|+{\mathcal M}^{\frac{m-1}2}|t-s|^\frac12\r)^\frac1{m-1}
\end{equation}
for every $(x,t), (y,s) \in \mathcal{K}_T \df = \mathcal K \times [0,T]$. 
\end{theorem}
\begin{remark}\label{Rmk:int:dist}
{\normalfont When dealing with the interior H\"older continuity of solutions to \eqref{Eq:1:5:3}, the statement is usually written as in Theorem~\ref{Thm:3:16:1}, namely in terms of the dimensionless, parabolic distance
\begin{equation}\label{def:par:dist}
d((x_1,t_1),(x_2,t_2))\df=\frac{|x_1-x_2|+{\mathcal M}^{\frac{m-1}2}|t_1-t_2|^{\frac12}}{\inf_{\ttop{(x,t)\in Q}{(y,s)\in\partial_p Q'}}
\Big(|x-y|+{\mathcal M}^{\frac{m-1}2}|t-s|^{\frac12}\Big)}.
\end{equation}
Unfortunately, this is not possible here, since typically the quantity at the denominator of \eqref{def:par:dist} is included in the constant $C_o$.
 A \emph{clean} statement could be given, if one were to accurately trace all the dependencies affecting $C_o$ in the arguments leading to the gradient bound.}
\end{remark}
\begin{proof}
Without loss of generality, we can assume ${\mathcal M}>1$. In order to simplify the exposition, we will study a rescaled solution, namely
\begin{equation*}
\tilde u(x,t) \df = \frac{u(x, {\mathcal M}^{1-m}t)}{{\mathcal M}}.
\end{equation*}
It is easy to check that this is, indeed, another solution of equation~\eqref{Eq:1:5:3}. We also begin with the auxiliary assumption that $\Omega_T$ is simply connected.

By the bounds in~\eqref{cont_bounds} we have 
\begin{equation}\label{cont_bounds2}
\tilde u(x,t) \le 1 \quad\text{and}\quad |\nabla \tilde u(x,t)^{m-1}| \le C_o,
\end{equation}
for every $(x,t) \in  \Omega_{\widetilde T} \df = \{(y,s) \in \mathbb R^{N+1} \mid (y, {\mathcal M}^{1-m}s) \in  \Omega_T\}$.
Due to the gradient bound here, we have that $\tilde u(\cdot, t)$ is locally $\frac1{m-1}$--H\"older continuous, uniformly in time, that is
\begin{equation}\label{spatial_holder}
\begin{split}
|\tilde u(x, t) - \tilde u(y, t)| &\le |\tilde u^{m-1}(x,t)-\tilde u^{m-1}(y,t)|^\frac1{m-1} \\
&\le  \sup_{(z,s) \in \Omega_{\widetilde T}}|\nabla \tilde u^{m-1}(z,s)|^{\frac1{m-1}}|x-y|^{\frac1{m-1}} \\
&\le C_o^\frac1{m-1}|x-y|^{\frac1{m-1}}
\end{split}
\end{equation}
for all $(x, t), (y,t) \in \Omega_{\widetilde T}$. 

Consider $(x_o, t_o)$ and $(x_o, t_1)$ two arbitrary points in $ \mathcal K_{\widetilde T}$, with $t_o<t_1$, where $$\mathcal K_{\widetilde T} \df= \mathcal K \times [0,\widetilde T] \ni (x_o, t_o). $$ 
It remains to show the H\"older condition in time, namely
\[
|\tilde u(x_o, t_o) - \tilde u(x_o, t_1)| \le \gamma_1(m,N,C_o)|t_o-t_1|^{\frac1{2(m-1)}}.
\]

By using the earlier notation ${\bf A}(\tilde u)=m\tilde u^{m-1}$, we define an operator $L$ by setting
\[
 L v \df= -\partial_t v+ {\bf A}(\tilde u)\Delta v + \nabla {\bf A}(\tilde u) \cdot \nabla v.
\]
Recall that the porous medium equation can now be written as
\[
L \tilde u =-\partial_t \tilde u + \Delta \tilde u^m =0.
\]
Notice also that in $\Omega_{\widetilde T}$ we have
\[
0 \le {\bf A}(\tilde u) \le m,
\] 
and moreover, we also know by~\eqref{cont_bounds2} that 
\[
|\nabla {\bf A}( \tilde u)| \le mC_o.
\]
Let $d_1\df=\dist(\mathcal K, \partial \Omega)$ be as defined earlier, and denote by $\mu_o=\mu_o(m,C_o)$ the maximum of the above upper bounds,
\[
\mu_o \df= \max\{m, mC_o\}.
\]
We define the set
\[
\mathcal N \df = B_r(x_o) \times (t_o ,t_1],
\]
where the radius $r \in (0, d_1)$ is to be chosen later. Observe that for such an $r$ we always have $\mathcal N \subset \Omega_{\widetilde T}$. Let
\[
\kappa \df= \sup_{t_o \le t \le t_1} |\tilde u(x_o, t)-\tilde u(x_o, t_o)|.
\] 
We need to show that $\kappa$ is bounded by $K|t_o-t_1|^\frac1{2(m-1)}$ where $K$ is a constant which depends only on the data $\{m,N\}$, $C_o$, and $d_1$, and in particular, is independent of $x_o$. 
We introduce the polynomial comparison functions $v^{\pm}:{\mathcal N}\cup\partial_p{\mathcal N}\to\rr$ defined by
\begin{align*}
v^{\pm}(x,t) \df= &\mu_o\l[1+\frac{2\kappa}{r^2}(N+r)\r](t-t_o) + \frac{\kappa}{r^2}|x-x_o|^2\\
&+ [C_or]^\frac1{m-1}\pm [\tilde u(x,t)-\tilde u(x_o, t_o)],
\end{align*}
and show that these functions are super-solutions for the operator $L$ in $\mathcal N$, while they attain non-negative values on the parabolic boundary $\partial_p \mathcal N$.  
By definition, in $\mathcal N$ we have
\[
Lv^\pm = -\mu_o-\frac{2\kappa}{r^2}(N\mu_o+\mu_or)+\frac{2\kappa}{r^2}\l[N{\bf A}(\tilde u)+\nabla {\bf A}(\tilde u) \cdot (x -x_o)\r] \le 0.
\]
Next, we study $v^\pm$ on the parabolic boundary of $\mathcal N$. For $t=t_o$ and $|x-x_o| < r$, by the H\"older estimate~\eqref{spatial_holder} we have that
\begin{align*}
v^\pm(x,t_o)=&\frac{\kappa}{r^2}|x-x_o|^2+[C_or]^\frac1{m-1}\pm [\tilde u(x,t_o)-\tilde u(x_o, t_o)] \ge 0,
\end{align*}
and for $|x-x_o|=r$ and $t_o < t \le t_1$ we obtain
\begin{align*}
v^\pm(x,t)&= \mu_o\l[1+\frac{2\kappa}{r^2}(N+r)\r](t-t_o) + \kappa+ [C_or]^\frac1{m-1}\pm [\tilde u(x,t)-\tilde u(x_o, t_o)] \\&\ge [C_or]^\frac1{m-1}\pm[\tilde u(x,t)-\tilde u(x_o, t)]+\kappa \pm[\tilde u(x_o,t)-\tilde u(x_o, t_o)] \\
&\ge 0,
\end{align*}
where in the last step, we used the definition of $\kappa$. 

Therefore, by the maximum principle for linear parabolic equations 
(see, for example, \cite{friedman-mono}, page 34), we deduce that $v^\pm \ge 0$ in 
$\overline {\mathcal N}$.  In particular, this implies $v^\pm(x_o, t) \ge 0$; in other words,
\[
\mp [\tilde u(x_o,t)-\tilde u(x_o, t_o)] \le [C_or]^\frac1{m-1}+\mu_o\l[1+\frac{2\kappa}{r^2}(N+r)\r](t-t_o).
\]
Taking the supremum over $t \in [t_o, t_1]$ yields
\begin{equation}\label{kappa_bound}
\kappa \le [C_or]^\frac1{m-1}+\mu_o(t_1-t_o)+\frac12\kappa\l[\frac{4\mu_o}{r^2}(N+r)(t_1-t_o)\r].
\end{equation}
We choose $r_o=r$ to be the positive root of 
\[
\frac{4\mu_o}{r^2}(N+r)(t_1-t_o) =1,
\]
that is
\[
r_o=2\mu_o(t_1-t_o)+2\l[N\mu_o(t_1-t_o)+\mu_o^2(t_1-t_o)^2\r]^\frac12.
\]
We need to have $r_o \in (0, d_1)$, which is true if 
\[
t_1-t_o < \frac{d_1^2}{4\mu_o(N+d_1)} \df=\delta_o.
\]
Inserting these in~\eqref{kappa_bound} yields
\begin{align*}
 \kappa &\le 2C_o^\frac1{m-1}\l[2\mu_o(t_1-t_o)+2\l[N\mu_o(t_1-t_o)+\mu_o^2(t_1-t_o)^2\r]^\frac12\r]^{\frac1{m-1}}\\
 &+2\mu_o(t_1-t_o) \\
&\le 2\l[C_o^\frac1{m-1}\big(2\mu_o\delta_o^{\frac12}+2[N\mu_o+\mu_o^2\delta_o]^{\frac12}\big)^{\frac1{m-1}}+\mu_o\delta_o^{\frac{2m-3}{2(m-1)}}\r](t_1-t_o)^\frac1{2(m-1)} \\
&\df=K(m,N,C_o, d_1)|t_1-t_o|^\frac{1}{2(m-1)}.
\end{align*}
Therefore, for every $t_o, t_1$, with $|t_o- t_1| < \delta_o$, we obtain that
\[
|\tilde u(x_o, t_1)-\tilde u(x_o, t_o)|\le K(m,N,C_o, d_1)|t_o-t_1|^{\frac1{2(m-1)}}.
\]
Since the constant $K$ does not depend on $x_o$, and $\mathcal K_{\widetilde T}$ is a compact set in time, we may now cover the set with a finite number of $t_k, k \ge 0$ for which $|t_k-t_{k-1}| < \delta_o$, to conclude that there exists a constant $\gamma_2=\gamma_2(m,N,C_o, d_1)$ such that for every $(x,t), (x,s) \in \mathcal K_{\widetilde T}$ we have 
\[
|\tilde u(x, t) - \tilde u(x, s)| \le \gamma_2|t-s|^{\frac1{2(m-1)}}.
\]
Rescaling back to $u$, and combining this with~\eqref{spatial_holder}, gives
\[
|u(x, t) - u(y, s)| \le \gamma {\mathcal M}\left(|x-y| + {\mathcal M}^{\frac{m-1}2}|t-s|^\frac12\right)^{\frac1{m-1}},
\]
for every $(x,t), (y,s) \in \mathcal K_T$ and for a constant $\gamma=\gamma(m, N, C_o, d_1)$. 

Finally, if $\Omega$ is not simply-connected, the result still follows by a standard covering argument for every compact $\mathcal K \Subset \Omega$.
\end{proof}

\section{The Proofs of Theorem~\ref{thm:main} and Corollary~\ref{Cor:hoelder} Concluded}\label{S:limit}
We will now conclude the proofs of our main results. First of all, notice that 
Proposition~\ref{prop:stability} is a {direct} consequence of Proposition~\ref{Ivanov_limit} below.  
Therefore, we limit ourselves to it.

The argument is quite standard: we will only provide a sketch of the proof, giving  detailed references for the missing arguments. We employ the following proposition due to Ivanov~\cite{Ivan95}.

\begin{proposition}\label{Ivanov_limit}
Let $\alpha \in (0,1)$ be the H\"older exponent of $u$ and suppose $u_\eps$ is the classical solution to problem~\eqref{regularised_equation} in $Q^*$. Then, there exists a subsequence $(\eps_k)_{k=1}^\infty$ such that
\[
u_{\eps_k} \to u \quad \text{in} \ \ C^\alpha(\overline {Q^*})
\]
and
\[
\nabla u_{\eps_k}^{\frac{m+1}2} \rightharpoonup \nabla u^{\frac{m+1}2}Ê\quad \text{weakly in}\ \ L^2(Q^*).
\]
\end{proposition}

\begin{proof}
Let $Q^*\df = B^* \times I^*$. Since $u \in C^\al({\overline{Q^*}})$ and $u^{\frac{m+1}2} \in L^2(I^*; W^{1,2}(B^*))$, by the setup of the regularised problem we have that
\[
u_\eps, u \in C^\alpha(\overline{Q^*}),\quad u_\eps^{\frac{m+1}2}, u^{\frac{m+1}2} \in L^2(I^*; W^{1,2}(B^*)),
\]
as well as
\[
0 \le u_\eps, u \le M,
\]
for a constant $M> 0$, uniformly for all $\eps>0$. By the Arzel\`a-Ascoli theorem, there exists a uniformly convergent subsequence $u_{\eps_k} \to v_1$ in $C^\alpha(\overline {Q^*})$ for some $v_1 \in C^\alpha(\overline {Q^*})$. On the other hand, the uniform boundedness of the functions $u_\eps^\frac{m+1}2$ in $L^2(I^*; W^{1,2}(B^*))$ implies that there exists another weakly convergent subsequence for which 
\[
u_{\eps_{k_i}}^{\frac{m+1}2}\rightharpoonup v_2^{\frac{m+1}2}\quad\text{weakly in}\ \  L^2(I^*; W^{1,2}(B^*)).
\]
For simplicity we drop the subscript $k_i$ and denote this subsequence by $u_\eps$. 

Since the convergence in $C^\alpha$-norm is pointwise uniform, we necessarily have $v_1=v_2$. Due to the uniform bounds that we have in the above function spaces, a tedious but straightforward argument shows that $v$ is a weak solution of the Dirichlet problem in $Q^*$ with boundary values given by $u$. By the local uniqueness of the solution to the equation, a direct consequence of {Proposition~\ref{prop:comparison}}, we obtain $v=u$. This finishes our sketch of the proof. For a rigorous treatment with all the details we refer to~\cite[Theorem~6.1 and Theorem~7.1]{Ivan95}.
\end{proof}

\subsection{Proof of the Main Results}
We can finally prove our main results concerning the boundedness of $|\nabla u^{m-1}|$, the consequent sharp H\"older estimate, and the instantaneous regularisation.

\begin{proof}[Proofs of Theorem~\ref{thm:main} and Corollary~\ref{Cor:hoelder}]

Since $u_\eps $ and $|\nabla u_\eps^{m-1}|$ are uniformly bounded, we obtain that for every $p>1$ there exists a function
$$
v_p \in L_{loc}^p(I^*, W_{loc}^{1,p}(B^*))
$$
such that 
\[
u_{\eps_k}^{m-1} \rightharpoonup v_p^{m-1}Ê\quad \text{weakly in}\ \ L_{loc}^p(I^*; W_{loc}^{1,p}(B^*)),
\]
where $u_{\eps_k}$ is a suitable subsequence of the sequence obtained in Proposition~\ref{Ivanov_limit}. By Proposition~\ref{Ivanov_limit} we have $u_{\eps_k} \to u$ in $C^\alpha(\overline {Q^*})$, which again guarantees that $v\df=v_p=u$ for all $p>1$. 

By the lower semicontinuity of the $L^p$--norm we obtain that for every $p >1$
\begin{equation}\label{grad_bound_finalized}
\begin{split}
\|\nabla u^{m-1}\|_{L^p} = \|\nabla v^{m-1}\|_{L^p}  &\le \liminf_{k\to 0} \|\nabla u_{\eps_k}^{m-1}\|_{L^p} \\
&\le \liminf_{k\to 0} \|\nabla u_{\eps_k}^{m-1}\|_{L^\infty} \le C(m,N, d_o, Q^*).
\end{split}
\end{equation}
Now $\|\nabla u^{m-1}\|_{L^p}$ is a bounded increasing sequence in $p$ and, thus, it has the limit
\[
\lim_{p\to \infty}\|\nabla u^{m-1}\|_{L^p}= \|\nabla  u^{m-1}\|_{L^\infty}.
\]
Taking the limit in~\eqref{grad_bound_finalized} gives the required gradient estimate for $u$.

Finally, the H\"older estimate of Lemma~\ref{holder} is uniform in $\eps$, and by the Arzel\`a-Ascoli theorem we have that 
\[
u_\eps \to v\quad \text{for some}\ \ vÊ\in C_{loc}^\frac1{m-1}(Q^*).
\] 
On the other hand, by Proposition~\ref{Ivanov_limit} we know that $u_\eps \to u$ in $C^\alpha(\overline {Q^*})$, which implies $v=u$. This finishes the proof.
\end{proof}

\section{Proof of Theorem~\ref{inst_reg_main}}\label{S:Inst}
Fix $\pto\in Q\Subset Q^*\Subset\Om_T$. Suppose that
\begin{equation}\label{eq:gamma:r}
{\dashint_{B_{r}(x_o)} u(x, t_o) \, dx \le \gamma r^{\frac1{m-1}}}
\end{equation}
for some constant $\gamma>0$, and for every $r$ small enough. {We now fix such an $r$}. Notice that, by the local continuity of $u$, we have $u\pto=0$.

We consider the $\eps$-regularisation \eqref{regularised_equation} of $u$ in $Q^*$
with {$\eps$ small enough, so that by Proposition~\ref{prop:stability} we have $\|u_\eps-u\|_{C^o(\overline {Q^*})}\le r^{\frac1{m-1}}$.}
{Moreover, let $u_{\eps, r} = \max\{u_\eps, r^{\frac1{m-1}}\}$, which is a subsolution, as it can be retrieved, for example, from Lemma~5.1 of \cite{DBGV-mono}.}

{We apply Lemma~\ref{Lm:sup-DiBe} with}
\begin{equation*}
\rho = r, \quad \sig=\frac12,\quad \text{and}\quad {\theta = \l[\sup_{t \in I^*}\dashint_{B_{r}(x_o)} u_{\eps,r}(x, t) \, dx\r]^{1-m}r^2}.
\end{equation*}
{Denoting
\begin{equation*}
\bar\theta \df =\l[\sup_{t \in I^*}\dashint_{B_{r}(x_o)} u_{\eps,r}(x, t) \, dx\r]^{1-m}
\end{equation*}
we obtain
\begin{equation}\label{eq:sup}
\begin{split}
&\sup_{\pto + Q^-_{\frac{r}2}({\bar\theta})} u_{\eps,r}(x, t)\\
&\le\kappa\l[\sup_{t \in I^*}\dashint_{B_{r}(x_o)} u_{\eps,r}(x, t) \, dx\r]^{1-m}
\left[\sup_{t_o-\theta<t<t_o}\dashint_{B_\rho(x_o)}u_{\eps, r}(x,t)\,dx\right]^m\\
&\qquad+\kappa\sup_{t \in I^*}\dashint_{B_{r}(x_o)} u_{\eps,r}(x, t) \, dx \\
&\le \kappa \sup_{t \in I^*}\dashint_{B_{r}(x_o)} u_{\eps,r}(x, t) \, dx,
\end{split}
\end{equation}
for a constant $\kappa$ depending only on the data $\{m,N\}$. 

Here we used the comparison principle to obtain
\begin{equation*}\label{rcomp}
 \sup_{t \in I^*}\dashint_{B_{r}(x_o)} u_{\eps, r}(x, t) \, dx \ge r^{\frac1{m-1}},
\end{equation*}
which implies
\begin{equation*}
\bar\theta r^2 =\theta\le r \le \frac {d_o}2,
\end{equation*}
and, thus,
\begin{equation*}\label{inclusion}
\pto + Q_{r}(\bar\theta) \subset Q^*.
\end{equation*}

In order to control the right-hand side of~\eqref{eq:sup} we estimate
\begin{equation*}\label{lower}
\begin{split}
\dashint_{B_{r}(x_o)} u_{\eps, r}(x, t) \, dx &\le  \dashint_{B_{r}(x_o)} u_{\eps}(x, t)+ r^{\frac1{m-1}}\, dx  \\
&\le  \dashint_{B_{r}(x_o)} u(x, t)+|u_{\eps}(x, t)-u(x, t)|+ r^{\frac1{m-1}}\, dx  \\
&\le  \dashint_{B_{r}(x_o)} u(x, t) \, dx  +2r^{\frac1{m-1}}\\
\end{split}
\end{equation*}
for every $t \in I^*$. By the Harnack inequality for the solution $u$, we conclude
\begin{equation*}
\sup_{t \in I^*} \dashint_{B_{r}(x_o)} u_{\eps, r}(x, t) \, dx \le \gamma_o  \dashint_{B_{r}(x_o)} u(x, t_o) \, dx+2r^{\frac1{m-1}}.
\end{equation*}

Plugging this into~\eqref{eq:sup} gives
\begin{align*}
\sup_{\pto + Q^-_{\frac{r}2}({\bar\theta})} u_{\eps,r}(x, t)
&\le \kappa  \dashint_{B_{r}(x_o)} u(x, t_o) \, dx+2\kappa r^{\frac1{m-1}},
\end{align*}
for a possibly modified constant $\kappa$ still depending only on the data $\{m,N\}$.

If $M$ is the quantity defined in \eqref{sigma_o-bis},
we conclude for all $r<1$ that $\theta \ge M^{1-m}r^2$. Therefore, using the comparison principle to deduce $u \le u_{\eps,r}$, and the assumption~\eqref{eq:gamma:r}, gives
\begin{equation*}
\sup_{\pto + Q^-_{\frac r2}(M^{1-m})} u(x,t) \le \kappa(\gamma+2)r^{\frac1{m-1}}.
\end{equation*}
%
%
%
%
%
%
This yields~\eqref{m-1_holder-main}, for the $r$ fixed in the beginning. Since $r$ was chosen arbitrarily, apart from being small, we conclude~\eqref{m-1_holder-main} for all $r$ small enough.}

Suppose now that there exists a sequence $r_{k_n} \to 0$ such that
\begin{align*}
\forall\,n\in\nn\qquad\dashint_{B_{r_{k_n}}(x_o)} u(x, t_o) \, dx \ge {\gamma}r_{k_n}^{\frac1{m-1}}.
\end{align*}
This implies
\begin{equation*}
\delta_n \df = \theta_{k_n}r_{k_n}^2 \le {\gamma^{1-m}}r_{k_n},\quad\theta_{k_n}\df=\left[\dashint_{B_{r_{k_n}}(x_o)} u(x, t_o) \, dx\right]^{1-m},
\end{equation*}
which gives us the sequence $\delta_n \to 0$. By Theorem~\ref{thm_quant}, we obtain the quantitative estimate in~\eqref{eq_quan-main}. Notice that {for $r_{k_n} \le1$} we have
\begin{align*}
C_o\theta_1^{\bar\al}\left[\frac{d_o}r\right]^{2}&=C_o d_o^{2}\frac{\theta_{k_n}^{\frac2{\al_o(m-1)}+1}r_{k_n}^{2\big[\frac2{\al_o(m-1)}+1\big]}}{r_{k_n}^{2}r_{k_n}^{2\big[\frac2{\al_o(m-1)}+1\big]}}\\
&\le C(m,N,d_o)\frac{\dl_n^{\frac2{\al_o(m-1)}+1}}{\gamma^{2\big[\frac1{\al_o}+m-1\big]}\dl_n^{4\big[\frac2{\al_o(m-1)}+1\big]}}=\frac {C(m,N,d_o)}{\gamma^{2\big[\frac1{\al_o}+m-1\big]}\dl_n^{3\bar\al}}.
\end{align*}
\section{Generalizations}\label{S:generalizations}

{The argument we presented relies on the intrinsic Harnack and weak Harnack inequalities, on the H\"older continuity of solutions in the interior, on a quantitative $L^\infty$ estimate, and on Lemma~\ref{first_ivanov_lemma}, but otherwise it is purely real analytic -- apart from passing to the limit $\eps \to 0$ in \S~\ref{S:limit}. The H\"older continuity of Lemma~\ref{hoelder}, the  Harnack inequalities of Theorems~\ref{Thm:3:15:1}--\ref{weak_harnack}, the quantitative $L^\infty$ estimate of Lemma~\ref{Lm:sup-DiBe},} hold also for doubly nonlinear equations of the form
\begin{equation}\label{general_operator}
u_t- \dvg \mathcal A(x,t, u, \nabla u)=0 \quad \text{weakly in $\Omega_T$,}
\end{equation}
where 
\begin{equation}\label{structure1}
\mathcal A(x,t, u, \nabla u)=m |u|^{m-1}|\nabla u|^{p-2}\nabla u,\quad m\ge2,\ p\ge2 
\end{equation}
or 
\begin{equation}\label{structure2}
\mathcal A_i(x,t, u, \nabla u)=m |u|^{m-1}|u_{x_i}|^{p-2} u_{x_1},\quad m\ge2,\ p\ge2. 
\end{equation}
For the precise statements and proofs of these result, we refer to~\cite{porzio-vespri, Ivanov-Mkrtychyan:1992,fornaro-sosio,vespri94}. 
The Barenblatt fundamental solution for equation \eqref{general_operator}--\eqref{structure1}
\begin{equation}\label{Barenblatt-dnl}
	{\mathcal B_{m,p}}(x,t)
	\df=
	\left\{
	\begin{array}{cl}
	\frac1{t^{\frac{N}{\lm}}}
	\left[
	1-b\left(\frac{|x|}{t^{\frac1\lm}}\right)^{\frac p{p-1}}
	\right]_+^{\frac{p-1}{m+p-3}},&t>0\\[10pt]
	0 &t\le 0
	\end{array}
	\right.
\end{equation}
where 
\begin{equation}\label{def:lambda-b-dnl}
	\lm=N(m+p-3)+p
	\quad
	\mbox{and}
	\quad
	 b=b(m,N,p)=\frac{p-1}p\frac{m+p-3}{(m+p-2)\lm^{\frac1{p-1}}}\, ,
\end{equation} 
suggests that in this case the optimal H\"older continuity exponent should be
\[
\al=\frac{p-1}{m+p-3}.
\]
An approach like the one we developed here, in principle should allow to check, whether this is indeed the case or not. However, we will not pursue this issue any further here.

\end{document}